# Assessing fluctuations of long-memory environmental variables based on the robustified dynamic Orlicz risk

Hidekazu Yoshioka[1, *], Yumi Yoshioka[2]

[1] Japan Advanced Institute of Science and Technology, 1-1 Asahidai, Nomi 923-1292, Japan

[2] Shimane University, Nishikawatsu-cho 1060, Matsue 690-8504, Japan

[*] Corresponding author: yoshih@jaist.ac.jp

ORCID:0000-0002-5293-3246 (Hidekazu Yoshioka)

**Abstract**

Environmental variables that fluctuate randomly and dynamically over time, such as water quality indices, are considered to be stochastic. They exhibit sub-exponential memory structures that should be accounted for in their modeling and analysis. Furthermore, risk assessments based on these environmental variables should consider limited data availability, which may introduce errors, e.g., model misspecifications, into their modeling. In this study, we present a pair of risk measures to determine the exponential disutility of a generic environmental variable both from below and above. The generic environmental variable is modelled as an infinite-dimensional nonlinear as well as affine stochastic differential equation and its moments and sub-exponential autocorrelations are estimated analytically. Novel risk measures, called dynamic robustified Orlicz risks, are formulated subsequently, and long, sub-exponential memory is efficiently addressed using them. The worst-case upper and lower bounds of the disutility are identified in closed form from the Hamilton–Jacobi–Bellman equations associated with the Orlicz risks. Finally, the proposed methodology is applied to weekly water quality data in a river environment in Japan.

**Keywords**

water quality; long-memory processes; stochastic differential equations; model uncertainty; dynamic Orlicz risks; Hamilton–Jacobi–Bellman equation

**Funding statement:** This study was supported by the Japan Society for the Promotion of Science (grant number 22K14441).

**Data availability statement**: The data will be made available upon reasonable request by the corresponding author.

**Conflict of interest disclosure**: The authors have no conflict of interests.

**Declaration of generative AI in scientific writing**: The authors did not use generative AI technology to write this manuscript.

**Permission to reproduce material from other sources**: N.A.

## 1. Introduction

### 1.1 Problem background

River systems play various roles in sustaining aquatic environments, ecosystems, and human life. The environmental conditions of rivers affect the habitat quality of aquatic species living in them significantly [1]. As a result, improving degraded environments and their valuation have become urgent social issues in the modern world [2]. Mathematical models have played a vital role in the modeling and analysis of both water quality and quantity in river systems across the world [3-5].

Environmental variables, such as water quality indices and river discharge, fluctuate randomly and dynamically over time. These variables are, therefore, described as stochastic processes; examples include, but are not limited to, discharge [6,7], dissolved oxygen [8], sulfate [9], and microbial communities [10,11]. In particular, stochastic differential equations (SDEs) are fundamental mathematical tools used for the description, analysis, and optimization of stochastic dynamical systems in general [12].

In environmental studies, SDEs have been employed to reproduce the time series of variables such as river discharge [13,14], sulfur dioxide concentration [15], total amount of carbon sequenced [16], antibiotic resistance [17], and generic pollutant indices that vary randomly with respect to time [18]. These studies utilize low-dimensional linear or affine SDEs whose autocorrelation functions decay exponentially, whereas the real-time series data of environmental variables exhibit sub-exponentials whose long memories decay only at algebraic speeds [19-22]. Similar phenomena have been reported for groundwater level fluctuations [23] and air pollution levels [24].

Stochastic processes with sub-exponential memories are often called long-memory processes or non-Markovian processes and have been represented using infinite-dimensional SDEs, such as lifted stochastic Volterra models [25,26] and superposed processes [27,28]. The two approaches use different noise structures—the former is based on fractional convolution of a single noise trajectory in time, while the latter is based on superposition (i.e., integration in a phase space) of infinitely many SDEs with different time scales. However, both face a common issue stemming from the infinite dimensionality of the resulting SDE. This difficulty is overcome by truncating the convolution or superposition at a finite degree of freedom to obtain a finite-dimensional SDE that approximates the original infinite-dimensional one [29,30]. This procedure is often called Markovian lifts or Markovian embedding, and has been applied to jump-driven SDE [19,31]. However, in the context of environmental variables, a unified theory to cover diffusion, jump, and jump-diffusion SDEs remains unavailable.

Modeling and analysis of environmental variables suffer from the issue of model uncertainty, e.g., modeling errors and model misspecification, because of limitations in data quality and/or quantity [32-36]. In the context of stochastic models, model uncertainty can be considered to be a distortion of the probability measure between the nominal model (i.e., the model identified based on the data) and the true model [37]. The relative entropy (or divergence) between the corresponding probability measures was proposed as a tool to measure the difference between the nominal and true models [38,39]. As environmental variables fluctuate significantly over time, a statistical index that can effectively capture their dynamic nature is required to address the issue of model uncertainty. However, such approaches are uncommon for environmental variables, particularly those with long memories.

Recently, the robustified dynamic Orlicz premium [40] has been proposed in insurance and financial research fields as a unified approach for the evaluation of the risk (i.e., the occurrence of extremely large or small premia) of stochastic processes under model uncertainty. Although this notion was not originally defined to assess the stochastic processes of environmental variables, it has recently been applied to engineering problems such as aquatic vegetation management [41] and photovoltaic power generation [42]. An advantage of the robustified dynamic Orlicz premium, henceforth called dynamic Orlicz for sake of simplicity, is its time consistency resembling that of dynamic programming [43] in a nonlinear manner. The risk under model uncertainty can be evaluated based on this by solving an associated optimality equation called the Hamilton–Jacobi–Bellman (HJB) equation. Formally, Orlicz risk is therefore a dynamic programming principle based on a nonlinear expectation called recursive utility, as utilized in economics and related research fields [44-47]. However, applications of recursive utilities, including the dynamic Orlicz risk have not yet been studied against environmental variables such as water quality indices.

### 1.2 Objective and contributions

The objective of this study is to apply the Orlicz risk to a jump-diffusion SDE model of a generic environmental variable with long memory. The target SDE is based on the superposition formalism [27] containing nonlinear as well as affine stochastic processes of the Cox–Ingersoll–Ross (CIR) type with jumps [48,49]. The use of the CIR-type SDE was motivated from the nonnegativity of its solutions that will be relevant for representing environmental variables such as the concentration of some pollutant that may accidentally increase following a jump process, and the affine nature with which the moment-generating functions and solutions to HJB equation can be obtained in closed forms under the novel superposition employed in this study. This jump-diffusion SDE is sufficiently flexible to represent the temporal evolution of water environmental variables because it represents the superposition of the Ornstein–Uhlenbeck processes [50] for jump-driven river discharge as well as diffusive processes as a special case as demonstrated later. The notion of superstatistics [51,52], in which parameters of the probability density depend on another random variable, can be considered to be a static counterpart of Markovian lifts. The superposition formalism enables the efficient capture of subexponential memories of water environmental variables. In addition, it can avoid computationally demanding statistical simulation of long-memory processes thanks to the use of a tractable HJB formalism.

In this study, the Orlicz risks provide both upper and lower bounds for an exponential (dis)utility, or the moment-generating function, of generic environmental variable, serving as a risk index for the target random variable owing to its convexity [53-55]. We demonstrate that the proposed SDE harmonizes with Orlicz risks, such that the resulting HJB equations admit closed-form solutions that can be implemented in engineering applications without resorting to sophisticated numerical methods. The key is the use of Markovian lifts to compute Orlicz risks via tractable finite-dimensional truncations of the HJB equations. A sufficient condition to well-pose Orlicz risks, depending on the functional forms of risk and uncertainty aversion, is also presented. Finally, the proposed mathematical framework is applied to real-world weekly water quality data in a Japanese river environment. Thus, this study consistently addresses the formulation, analysis, and application of the new SDE and recursive utility in the field of environmental research.

The remainder of this paper is organized as follows. The infinite-dimensional SDE analyzed in this study is explained in **Section 2**. The Orlicz risks are presented in **Section 3** and the corresponding HJB equations and closed-form solutions are derived. The application of the proposed mathematical framework to water environmental data of a Japanese river is described in **Section 4**. Finally, the paper is concluded in **Section 5** and directions of future research are discussed. **Appendices** contain the proofs of the propositions stated in the main text and several auxiliary results.

## 2. Stochastic model

The SDE used in this study is described in this section. We first discuss the one-dimensional counterpart, which is simpler and more tractable, while sharing certain properties with the infinite-dimensional one.

### 2.1 One-dimensional model

We work in a complete probability space $(\Omega, \mathcal{F}, \mathbb{P})$ with $\Omega$ being the set of all possible events, the sigma algebra $\mathcal{F}$ generated by jump and continuous noises (Brownian, pure-jump Lévy processes, and their multi-dimensional versions that will appear later), and the associated probability measure $\mathbb{P}$, as per the conventional studies [12]. The 1-D standard Brownian motion at time $t \geq 0$ is denoted by $B_t$. The pure-jump Lévy process with bounded-variation jumps, also called a subordinator, with the Lévy measure $v(\mathrm{d}z)$ is denoted by $L_t$ at $t \geq 0$. This $L$ is independent of $B$. The processes $B$ and $L$ play different roles with each other in our model although both of them are Lévy processes. This is the reason that we deal with them separately. The Lévy measure $v$ is assumed to satisfy

$$M_k = \int_0^{+\infty} z^k v(\mathrm{d}z) < +\infty, \quad k \in \mathbb{N}. \tag{1}$$

The quantities $M_K$ ($k \in \mathbb{N}$) represent moments of the jump size up to a proportional constant. The boundedness assumption gives a sufficient condition to guarantee the existence of the moments of our SDE (e.g., (45)-(47)). The condition (1) is not restrictive because the typical Lévy measure $v(\mathrm{d}z)$ used in applications exhibits the tempered stable form [19,31,50]:

$$v(\mathrm{d}z) = \frac{\gamma}{z^{1+\alpha}} \exp(-\beta z)\mathrm{d}z, \quad z > 0 \tag{2}$$

with parameters $\beta, \gamma > 0$ and $\alpha < 1$. This $v$ covers compound Poisson cases with exponential or gamma density ($\alpha < 0$) and subordinates with infinite activity ($\alpha \geq 0$) [56,57]. Our setting is therefore sufficiently general to deal with these specific cases. This is an advantage of our mathematical framework.

The one-dimensional Itô's SDE that governs the continuous-time process $X = (X_t)_{t \geq 0}$ representing an environmental variable, such as the concentration of a pollutant, is set as (e.g., [48,49])

$$\mathrm{d}X_t = (a - rX_t)\mathrm{d}t + \sigma\sqrt{rX_t}\mathrm{d}B_t + \mathrm{d}L_t,\ \ t > 0 \tag{3}$$

with an initial condition $X_0 \geq 0$ and parameters $a, r, \sigma > 0$. This SDE, called the jump CIR (JCIR) process, admits a unique pathwise continuous solution that is nonnegative and has the following exponential autocorrelation function

$$\mathrm{ACF}(h) = \exp(-rh) \tag{4}$$

for time lag $h \geq 0$. Therefore, the parameter $r$ serves as a reversion speed towards the mean.

The stationary moment-generating function of $X$ is obtained in a closed form as follows:

$$\mathbb{E}\left[\exp(pX_t)\right] = \exp\left(\int_0^{+\infty}\left\{au_r(s) + \int_0^{+\infty}\left(\exp(u_r(s)z) - 1\right)v(\mathrm{dz})\right\}\mathrm{ds}\right),\ \ p < \min\left\{\frac{2}{\sigma^2}, \beta\right\}, \tag{5}$$

where

$$u_r(s) = \left(\frac{\sigma^2}{2} + \left(\frac{1}{p} - \frac{\sigma^2}{2}\right)e^{rs}\right)^{-1},\ \ s > 0. \tag{6}$$

The moments of $X$ can also be obtained analytically through $(5)$. The infinite-dimensional SDE introduced in the next subsection inherits this remarkable property, which is leveraged in **Section 4**. Actually, $(5)$ is a special case of **Proposition 1** presented later ($n = 1$, $\Phi(x) = x$, $\lambda_\phi, \lambda_\varphi \to +0$, $t \to -\infty$ there) whose detailed derivation procedure is found in **Appendix A**. The quantity $p$ in $(5)$ is the real parameter to define the moment-generating function. Each moment can be formally obtained by considering the derivative at $p \to 0$ (e.g., $\mathbb{E}[X_t] = \left(\frac{\partial}{\partial t}\mathbb{E}\left[\exp(pX_t)\right]\right)_{p \to 0}$).

## 2.2 The superposed model

The infinite-dimensional SDE generalizing SDE $(3)$ is formulated using the superposition formalism [27]. The SDE, referred to as the superposition of the JCIR processes (supJCIR process), is essentially the sum of infinitely many mutually independent JCIR processes with mutually different reversion speeds, $r$. The coexistence of processes with mutually different $r$ emerges as a sub-exponential autocorrelation, as demonstrated below.

We assume that there exists a probability measure $\pi(\mathrm{d}r)$ that generates the reversion speed $r > 0$ and approximate it using a discrete probability measure.

$$\pi_n(\mathrm{d}r) = \sum_{i=1}^{n} c_i \delta(r - r_i) \tag{7}$$

with the Dirac's delta $\delta$; a positive, strictly increasing, and bounded sequence $\{r_i\}_{i=1,2,3,\ldots,n}$; and another positive sequence $\{c_i\}_{i=1,2,3,\ldots,n}$ with $\sum_{i=1}^{n} c_i = 1$. We assume that the sequences, $\{r_i\}_{i=1,2,3,\ldots,n}$ and $\{c_i\}_{i=1,2,3,\ldots,n}$, are selected such that $\pi_n$ converges to $\pi$ in the sense of distributions [58].

We construct the supJCIR process from its finite-dimensional counterpart as a finite sum of mutually independent JCIR processes. We set the degree-of-freedom $n \in \mathbb{N}$, and set $n$ mutually independent copies $B^{(i)}$ ($i = 1,2,3,...,n$) of the Brownian motion $B$ and $n$ mutually independent scaled copies $L^{(i)}$ ($i = 1,2,3,...,n$) of the Lévy process $L$, where the Lévy measure of $L^{(i)}$ is $c_i v$. $L^{(i)}$ is independent of $B^{(i)}$. The finite-dimensional version of the supJCIR process is referred to as $Y^{(n)}$. This variable represents the continuous-time evolution of an environmental variable and is given by

$$Y_t^{(n)} = \sum_{i=1}^{n} X_t^{(i)}, \quad t > 0 \tag{8}$$

with the Itô's SDE

$$\mathrm{d}X_t^{(i)} = \left(ac_i - r_i X_t^{(i)}\right)\mathrm{d}t + \sigma\sqrt{r_i X_t^{(i)}}\mathrm{d}B_t^{(i)} + \mathrm{d}L_t^{(i)}, \tag{9}$$

starting with the initial condition $X_0^{(i)}$ ($i = 1,2,3,...,n$). Representation (8) implies that the evolution of the environmental variable $Y^{(n)}$ is driven by several microscopic evolutions (9) at different characteristic speeds $r_i$. For major variables, e.g., the concentration of some nutrients and flow discharge, the existence of microscopic processes with different timescales can be attributed to slow and fast runoff processes in the watershed [59,60].

Based on the independence of each $X^{(i)}$, we obtain the stationary moment generating function

$$\begin{aligned}\mathbb{E}\left[\exp\left(pY_t^{(n)}\right)\right] &= \exp\left(\sum_{i=1}^{n}\int_0^{+\infty}\left\{ac_i u_r(s) + \int_0^{+\infty}\left(\exp\left(u_r(s)z\right)-1\right)c_i v(\mathrm{d}z)\right\}\mathrm{d}s\right) \\ &= \exp\left(R^{(n)}\int_0^{+\infty}\left\{au_1(s) + \int_0^{+\infty}\left(\exp\left(u_1(s)z\right)-1\right)v(\mathrm{d}z)\right\}\mathrm{d}s\right)\end{aligned}, \quad p < \min\left\{\frac{2}{\sigma^2}, \beta\right\} \tag{10}$$

with

$$R^{(n)} = \sum_{i=1}^{n}\frac{c_i}{r_i}. \tag{11}$$

The moment-generating function (10) of the finite-dimensional supJCIR process resembles that (5) of the JCIR process—their difference is the multiplier of the integration within the argument of "exp." Hence, the finite-dimensional version inherits the closed-form availability of the moment-generating function in the JCIR process. Further, the autocorrelation function associated with (8) is:

$$\mathrm{ACF}(h) = \frac{1}{R^{(n)}}\sum_{i=1}^{n}\frac{c_i}{r_i}\exp\left(-r_i h\right), \quad h \geq 0, \tag{12}$$

which is the weighted sum of the exponential equation (4), reflecting the existence of multiple timescales.

Let us assume that we have the following inverse moment condition:

$$R = \int_0^{+\infty}\frac{\pi(\mathrm{d}r)}{r} < +\infty. \tag{13}$$

This indicates that the probability measure $\pi$ is sufficiently regular near $r = 0$. Phenomenologically, this means that the proportion of JCIR processes with small $r$, which decays slowly, is sufficiently small. This assumption is not restrictive in applications where we use Gamma-type $\pi$ [50] in the following:

$$\pi(\mathrm{d}r) = g(r)\mathrm{d}r \sim r^{\omega-1}\exp\left(-\frac{r}{\theta}\right)\mathrm{d}r\,,\quad r>0 \tag{14}$$

with $\omega > 1$ and $\theta > 0$.

The form of (10) suggests that it has the limit

$$\begin{aligned}&\lim_{n\to+\infty}\mathbb{E}\left[\exp\left(pY_t^{(n)}\right)\right]\\&=\exp\left(R\int_0^{+\infty}\left\{au_1(s)+\int_0^{+\infty}\left(\exp\left(u_1(s)z\right)-1\right)v(\mathrm{d}z)\right\}\mathrm{d}s\right)\end{aligned},\quad p<\min\left\{\frac{2}{\sigma^2},\beta\right\}, \tag{15}$$

if $R^{(n)} \to R$ as $n \to +\infty$. This convergence is justified by choosing quantile-based $\{r_i\}_{i=1,2,3,\ldots,n}$ and $\{c_i\}_{i=1,2,3,\ldots,n}$ [19]. We define a stationary process, the supJCIR process $Y = (Y_t)_{t\geq 0}$ as a weak limit (limit in the sense of distributions) of the finite-dimensional process (8) under $n \to +\infty$, such that

$$\mathbb{E}\left[\exp(pY_t)\right]=\exp\left(R\int_0^{+\infty}\left\{au_1(s)+\int_0^{+\infty}\left(\exp\left(u_1(s)z\right)-1\right)v(\mathrm{d}z)\right\}\mathrm{d}s\right),\quad p<\min\left\{\frac{2}{\sigma^2},\beta\right\}. \tag{16}$$

The convergence $R^{(n)} \to R$ also suggests the convergence of the autocorrelation function:

$$\frac{1}{R^{(n)}}\sum_{i=1}^{n}\frac{c_i}{r_i}\exp(-r_i h)\to\frac{1}{R}\int_0^{+\infty}\frac{\exp(-rh)\pi(\mathrm{d}r)}{r},\quad h\geq 0. \tag{17}$$

For Gamma-type $\pi$ in particular, we obtain the following subexponential autocorrelation function of $Y$:

$$\mathrm{ACF}(h)=\frac{1}{(1+\theta h)^{\omega-1}},\quad h\geq 0. \tag{18}$$

## 3. Orlicz risks

### 3.1 Orlicz function and relative entropy

The Orlicz risk in this study provides the upper and lower bounds of the exponential disutility under model uncertainty. Hence, it can be used to assess the worst-case environmental variables that are overestimated or underestimated owing to model misspecification.

We first focus on the finite-dimensional model (9) of $Y^{(n)}$, and later infer the infinite-dimensional limit, because it deals directly with the infinite-dimensional SDE of $Y$ (see **Appendix C**). In this study, the modeling errors are taken to be distortions of the drift and jumps in (9). The modelling errors are evaluated using Radon–Nikodým derivatives by invoking the jump-diffusion structure of (9) following the methodology of robust control [37,61,62].

An Orlicz risk has two key elements—an Orlicz function and a relative entropy [40]. An Orlicz function is a nonnegative, continuous, strictly increasing, and convex function $\Phi:[0,+\infty)\to[0,+\infty)$, satisfying $\Phi(0)=0$ and $\Phi(1)=1$. Typical examples include power functions $\Phi(x)=x^m$ ($m>1$), exponential functions $\Phi(x)=\frac{e^{mx}-1}{e^{m}-1}$ ($m>0$), and their weighted sums. Few assumptions on the Orlicz function will be added in each proposition when necessary. The Orlicz function will be fully specified in our application.

Let us set a measurable vector process $\phi=(\phi_t)_{t\geq 0}=\left(\left\{\phi_t^{(i)}\right\}_{i=1,2,3,\ldots,n}\right)_{t\geq 0}$ and a distorted Brownian motion $W=\left(\left\{W_t^{(i)}\right\}_{i=1,2,3,\ldots,n}\right)_{t\geq 0}$ as

$$\mathrm{d}W_t^{(i)}=\mathrm{d}B_t^{(i)}-\phi_t^{(i)}\mathrm{d}t\,. \tag{19}$$

Let us consider a measurable vector field process $\varphi=\left(\left\{\varphi_t^{(i)}(\cdot)\right\}_{i=1,2,3,\ldots,n}\right)_{t\geq 0}$, with each $\varphi_t^{(i)}(\cdot)$ being a positive mapping. Then, consider the distorted Lévy process $J=\left(\left\{J_t^{(i)}\right\}_{i=1,2,3,\ldots,n}\right)_{t\geq 0}$ where $J_t^{(i)}$ is the subordinator $L_t^{(i)}$ whose Lévy measure is modulated as $c_i\varphi_t^{(i)}(z)v(\mathrm{d}z)$. The probability measure under which the processes $W$ and $J$ are martingales is denoted by $\mathbb{Q}(\phi,\varphi)$. The probability measure without any model uncertainty is denoted by $\mathbb{P}=\mathbb{Q}(\mathbf{0},\mathbf{1})$. Indeed, $W=B$ and $L=J$ when $(\phi,\varphi)$ are constant processes $(\mathbf{0},\mathbf{1})=\left(\{0\}_{i=1,2,3,\ldots,n},\{1\}_{i=1,2,3,\ldots,n}\right)$. The expectation under $\mathbb{Q}(\phi,\varphi)$ is denoted by $\mathbb{E}_{\mathbb{Q}(\phi,\varphi)}$.

Given some $(\phi,\varphi)$, the difference between the measures $\mathbb{P}$ and $\mathbb{Q}(\phi,\varphi)$ is measured in terms of relative entropies [37]. Concerning the diffusive part, the relative entropy (Kullback–Leibler divergence) during the time interval $(\underline{t},\overline{t})$ with $0\leq\underline{t}\leq\overline{t}$ is

$$R_{\phi,(\underline{t},\overline{t})}=\mathbb{E}_{\mathbb{Q}(\phi,\varphi)}\left[\frac{1}{2}\sum_{i=1}^{n}\int_{\underline{t}}^{\overline{t}}\left(\phi_s^{(i)}\right)^2\mathrm{d}s\right]. \tag{20}$$

Similarly, for the jump part, the generalized relative entropy (Tsallis divergence) with a parameter $q>0$ during the time interval $(\underline{t},\overline{t})$ with $0\leq\underline{t}\leq\overline{t}$ is

$$\begin{aligned}R_{\varphi,(\underline{t},\overline{t})}&=\mathbb{E}_{\mathbb{Q}(\phi,\varphi)}\left[\sum_{i=1}^{n}\int_{\underline{t}}^{\overline{t}}\int_{0}^{+\infty}\frac{1}{1-q}\left(1-q-\left(\varphi_s^{(i)}(z)\right)^q+q\varphi_s^{(i)}(z)\right)v(\mathrm{d}z)\mathrm{d}s\right]\ (q\neq 1)\\&=\mathbb{E}_{\mathbb{Q}(\phi,\varphi)}\left[\sum_{i=1}^{n}\int_{\underline{t}}^{\overline{t}}\int_{0}^{+\infty}\left(\varphi_s^{(i)}(z)\ln\varphi_s^{(i)}(z)-\varphi_s^{(i)}(z)+1\right)v(\mathrm{d}z)\mathrm{d}s\right]\ (q=1)\end{aligned}. \tag{21}$$

Tsallis divergence is used for the jump part as Kullback–Leibler divergence ($q=1$) fails in applications of environmental variables with (2) (See **Proposition 1**). For later use, given uncertainty aversion coefficients $\lambda_\phi,\lambda_\varphi>0$, we set the following weight sum of the relative entropies:

$$R_{(\underline{t},\overline{t})}=\frac{1}{\lambda_\phi}R_{\phi,(\underline{t},\overline{t})}+\frac{1}{\lambda_\varphi}R_{\varphi,(\underline{t},\overline{t})}\,. \tag{22}$$

Loosely speaking, choosing a larger uncertainty aversion coefficient $\lambda_\phi$ or $\lambda_\varphi$ implies larger model uncertainty measured by relative entropy.

### 3.2 Formulation for the upper bound

Recursive utilities are obtained as worst-case disutilities. First, we explain the Orlicz risk for the upper bound originally proposed by [40]. The information available at time $t \geq 0$ is denoted by $\mathcal{F}_t$ and is generated by $B_t$ and $L_t$. At time $t$, the Orlicz risk $\|Z\|_{\Phi,t}$ of a positive random variable $Z$, measurable with respect to $\mathcal{F}_t$, is defined as

$$\|Z\|_{\Phi,t} = \inf\left\{ h>0 \left| \sup_{\mathbb{Q}(\phi,\varphi)} \mathbb{E}_{\mathbb{Q}(\phi,\varphi)}\left[ \underbrace{\Phi\left(\frac{Z}{h}\right)}_{\text{Risk aversion}} - \underbrace{C(\phi,\varphi)}_{\text{Uncertainty aversion}} \middle| \mathcal{F}_t \right] \leq 1,\ h \text{ is } \mathcal{F}_t\text{-measurable} \right.\right\}, \tag{23}$$

where the measure of model uncertainty $C(\phi,\varphi) \geq 0$ is $\mathcal{F}_t$-measurable. Equation (23) is abstract, but has a probabilistic interpretation that $\|Z\|_{\Phi,t}$ is the worst-case upper bound of $Z$ subject to the model uncertainty $C(\phi,\varphi)$. Indeed, the maximization using "sup" (23) implies the worst-case situation and the second term in the expectation is interpreted as a penalization of the model uncertainty; uncertainty is evaluated more seriously for larger $C(\phi,\varphi)$. The convexity of the $\Phi$ suggests that the maximizing $h$ is estimated to be greater for larger $Z$. Indeed, the classical Jensen inequality implies that

$$\Phi^{(-1)}\left(\mathbb{E}_{\mathbb{Q}(\phi,\varphi)}\left[\Phi\left(\frac{Z}{h}\right)\middle|\mathcal{F}_t\right]\right) \geq \mathbb{E}_{\mathbb{Q}(\phi,\varphi)}\left[\frac{Z}{h}\middle|\mathcal{F}_t\right]. \tag{24}$$

This suggests that $\Phi$ inside expectations yields an overestimation of the random variable of interest. The second term in the expectation of (23) combined with the supremum suggests that the first term in the expectation is estimated to be larger than that without uncertainty, and hence represents uncertainty aversion. From an alternative perspective related to optimization problems, this term can also be considered to be a penalty against data overfitting.

We now adapt the Orlicz risk (23) to our model. Assuming that we want to evaluate an exponential disutility $\mathbb{E}_{\mathbb{Q}(\phi,\varphi)}\left[\exp\left(pY_T^{(n)}\right)\middle|\mathcal{F}_t\right]$ ( $p>0$ ) at time $T>0$ based on the information available at $t \leq T$; for example, one may evaluate the worst-case expected upper bound for environmental pollution. The time increment $k>0$ is set to be smaller than $T-t$. The worst-case upper bound of the exponential disutility estimate $\bar{\Psi}_t > 0$, considering model uncertainty, is inspired by (23):

$$\sup_{\mathbb{Q}(\phi,\varphi)}\left\{\mathbb{E}_{\mathbb{Q}(\phi,\varphi)}\left[\Phi\left(\frac{\bar{\Psi}_{t+k}}{\bar{\Psi}_t}\right)\middle|\mathcal{F}_t\right] - R_{(t,t+k)}\right\} = 1,\ \ t<T \tag{25}$$

with

$$\bar{\Psi}_T = \exp\left(pY_T^{(n)}\right). \tag{26}$$

Equation (25) is a recursion that can determine the worst-case upper bound $\bar{\Psi}$ backwards in time from the terminal value (26). The first term of (25) inside the expectation represents risk aversion. Due to the increasing nature and convexity of $\Phi$, larger values of $\bar{\Psi}_{t+k}$ correspond to larger values of $\bar{\Psi}_t$, which contribute to the overestimation of exponential disutility. $R_{(t,t+k)}$ was substituted into $C(\phi,\varphi)$ to account for model uncertainty using relative entropies.

For subsequent use, the first and second derivatives of $\Phi$ are denoted by $\Phi'$ and $\Phi''$, respectively. Let us set $\bar{\lambda} = \frac{\Phi'(1)^2 \lambda_\phi + \Phi''(1)}{\Phi'(1)} > 0$ and introduce the *q*-exponential functions [63]:

$$\exp_q(z) = \begin{cases} \left(1+(1-q)z\right)^{\frac{1}{1-q}} & \left(1+(1-q)z > 0 \text{ and } q \neq 1\right) \\ \exp(z) & (q=1) \end{cases}. \tag{27}$$

We also use the Kronecker delta $\Delta_{i,j}$ ($\Delta_{i,j} = 1$ if $i = j$, and $\Delta_{i,j} = 0$ otherwise).

We obtain the HJB equation for the upper bound $\bar{\Psi}$ and its solution, as given in **Proposition 1** (Proof is in **Appendix A**). This is the main theoretical result of this study. **Proposition 1** states that the proposed mathematical framework can be readily implemented if the benchmark model is identified and the Orlicz risk is specified. In this method, the resulting ordinary differential equations can be evaluated using common numerical methods, such as the forward Euler method. The integrals in these equations can be easily approximated using a common numerical method, such as trapezoidal or midpoint rule. **Appendix C** presents the technical details of the implementation.

***Proposition 1***

*Let us assume* $0 < q < 1$ *and* $0 < p < \min\left\{\frac{2}{\sigma^2\left(1+\bar{\lambda}\right)}, \beta\right\}$. *Moreover, let us also assume that the process* $\bar{\Psi}_t$ *has a Markovian form* $\bar{\Psi}_t = \bar{\Psi}\left(t, X_t^{(1)}, X_t^{(2)}, ..., X_t^{(n)}\right) = \bar{\Psi}\left(t, \left\{X_t^{(i)}\right\}_{i=1,2,3,...,n}\right)$ *with a sufficiently smooth* $\bar{\Psi}$. *In addition, assume that there exist constants* $C_0, C_1 > 0$ *and* $m \in \left(1, \frac{\beta}{p}(1-q)\right)$ *such that*

$$\Phi(x) \le C_0 x^m + C_1, \quad x \ge 0. \tag{28}$$

*Then, the HJB equation that governs* $\bar{\Psi} = \bar{\Psi}\left(t, \left\{x_i\right\}_{i=1,2,3,...,n}\right)$ *is given by*

$$\underbrace{\frac{\partial\bar{\Psi}}{\partial t}+\sum_{i=1}^{n}\left(\left(ac_i-r_ix_i\right)\frac{\partial\bar{\Psi}}{\partial x_i}+\frac{\sigma^2}{2}r_ix_i\frac{\partial^2\bar{\Psi}}{\partial x_i^2}\right)}_{\text{Benchmark model}}$$

$$+\underbrace{\frac{\bar{\lambda}}{\bar{\Psi}}\sum_{i=1}^{n}\frac{\sigma^2}{2}r_ix_i\left(\frac{\partial\bar{\Psi}}{\partial x_i}\right)^2}_{\text{Diffusion effect considering risk and uncertainty}} \qquad \text{for } \left(t,\{x_i\}_{i=1,2,3,\ldots,n}\right)\in(-\infty,T)\times\mathbb{R}^n \tag{29}$$

$$+\underbrace{\frac{1}{\lambda_\varphi\Phi'(1)}\sum_{i=1}^{n}c_i\int_0^{+\infty}\left\{\exp_q\left(\lambda_\varphi\left(\Phi\left(\hat{\Psi}_i\right)-1\right)\right)-1\right\}v(\mathrm{d}z_i)}_{\text{Jump effect considering risk and uncertainty}}=0$$

*with*

$$\hat{\Psi}_i=\frac{\bar{\Psi}\left(t,\left\{x_j+\Delta_{i,j}z_j\right\}_{j=1,2,3,\ldots,n}\right)}{\bar{\Psi}\left(t,\left\{x_j\right\}_{j=1,2,3,\ldots,n}\right)} \tag{30}$$

*subject to the terminal condition*

$$\bar{\Psi}\left(T,\{x_i\}_{i=1,2,3,\ldots,n}\right)=\exp\left(p\sum_{i=1}^{n}x_i\right). \tag{31}$$

*In addition, (29) admits a smooth solution*

$$\bar{\Psi}=\exp\left(\sum_{i=1}^{n}\bar{\rho}_t^{(i)}x_i+\bar{\tau}_t\right) \tag{32}$$

*with time-dependent coefficients* $\left(\bar{\tau}_t,\left\{\bar{\rho}_t^{(i)}\right\}_{i=1,2,3,\ldots,n}\right)_{t\le T}$ *such that*

$$\bar{\rho}_t^{(i)}=\frac{1}{\frac{\sigma^2\left(1+\bar{\lambda}\right)}{2}+\left(\frac{1}{p}-\frac{\sigma^2\left(1+\bar{\lambda}\right)}{2}\right)e^{r_i(T-t)}},\quad i=1,2,3,\ldots,n \tag{33}$$

*and*

$$\frac{\mathrm{d}\bar{\tau}_t}{\mathrm{d}t}+a\sum_{i=1}^{n}c_i\bar{\rho}_t^{(i)}+\frac{1}{\lambda_\varphi\Phi'(1)}\sum_{i=1}^{n}c_i\int_0^{+\infty}\left\{\exp_q\left(\lambda_\varphi\left(\Phi\left(e^{\bar{\rho}_t^{(i)}z_i}\right)-1\right)\right)-1\right\}v(\mathrm{d}z_i)=0,\quad t\le T. \tag{34}$$

***Remark 1*** We have excluded the case $q>1$ from the assumption of **Proposition 1** as the last term, i.e., the integral, in the HJB equation (29) diverges as $n\to+\infty$ if $q=1$. Therefore, Kullback–Leibler divergence cannot be used for the jump part under the infinite-dimensional limit, which is the primary case of interest in this study, even though it works in the finite-dimensional limit case. This is the main reason that we proposed to use the Tsallis divergence in this study so that the drawback of the Kullback–Leibler one can be completely avoided. See, also **Appendix A**.

***Remark 2*** **Proposition 1** is a theoretical result explicitly connecting superstatistics and the HJB equation.

### 3.3 Formulation of the lower bound

The lower bound of the exponential disutility is also computed, which provides an optimistic bound for disutility. Its formulation is a symmetric counterpart of that for the upper bound presented in the previous subsection; however, we assume in this subsection that $\Phi$ is strictly increasing as well as concave again with $\Phi(0)=0$ and $\Phi(1)=0$, so that the Orlicz risk can provide a lower bound for exponential disutility. A typical example will be $\Phi(x)=x^{1/m}$ with $m>1$.

For the lower bound, the minimization counterpart of (23) is symmetrically defined as:

$$\|Z\|_{\Phi,t}=\sup\left\{h>0\,\middle|\,\inf_{\mathbb{Q}(\phi,\varphi)}\mathbb{E}_{\mathbb{Q}(\phi,\varphi)}\left[\underbrace{\Phi\left(\frac{Z}{h}\right)}_{\text{Risk aversion}}+\underbrace{C(\phi,\varphi)}_{\text{Uncertainty aversion}}\middle|\mathcal{F}_t\right]\geq 1,\ h\text{ is }\mathcal{F}_t\text{-measurable}\right\}. \tag{35}$$

As in the upper bound case, we assume the worst-case lower bound of the exponential disutility estimate $\underline{\Psi}_t>0$, considering model uncertainty based on the equality:

$$\inf_{\mathbb{Q}(\phi,\varphi)}\left\{\mathbb{E}_{\mathbb{Q}(\phi,\varphi)}\left[\Phi\left(\frac{\underline{\Psi}_{t+k}}{\underline{\Psi}_t}\right)\middle|\mathcal{F}_t\right]+R_{(t,t+k)}\right\}=1,\quad t<T \tag{36}$$

with

$$\underline{\Psi}_T=\exp\left(pY_T^{(n)}\right). \tag{37}$$

Equation (36) is a backward recursion that can determine the worst-case lower bound $\underline{\Psi}$ based on the terminal value in (37) to an arbitrary time before $T$. For later use, we define $\underline{\lambda}=\frac{\Phi'(1)^2\lambda_\phi-\Phi''(1)}{\Phi'(1)}>0$.

We obtain the HJB equation for the lower bound and its solution, as given in **Proposition 2** (Proof is in **Appendix B**). This is also a theoretical result of this study.

***Proposition 2***

*Assume* $q\geq 1$ *and* $0<p<\min\left\{\frac{2}{\sigma^2},\beta\right\}$. *Assume that the process* $\underline{\Psi}_t$ *has a Markovian form* $\underline{\Psi}_t=\underline{\Psi}\left(t,X_t^{(1)},X_t^{(2)},\ldots,X_t^{(n)}\right)=\underline{\Psi}\left(t,\left\{X_t^{(i)}\right\}_{i=1,2,3,\ldots,n}\right)$ *with a sufficiently smooth function* $\underline{\Psi}$. *Then, the HJB equation that governs* $\underline{\Psi}=\underline{\Psi}\left(t,\{x_i\}_{i=1,2,3,\ldots,n}\right)$ *is given by*

$$\underbrace{\frac{\partial \underline{\Psi}}{\partial t}+\sum_{i=1}^{n}\left(\left(ac_i-r_ix_i\right)\frac{\partial \underline{\Psi}}{\partial x_i}+\frac{\sigma^2}{2}r_ix_i\frac{\partial^2 \underline{\Psi}}{\partial x_i^2}\right)}_{\text{Benchmark model}}$$

$$-\underbrace{\frac{\underline{\lambda}}{\underline{\Psi}}\sum_{i=1}^{n}\frac{\sigma^2}{2}r_ix_i\left(\frac{\partial \underline{\Psi}}{\partial x_i}\right)^2}_{\text{Diffusion effect considering risk and uncertainty}} \qquad \textit{for } \left(t,\{x_i\}_{i=1,2,3,\ldots,n}\right)\in(-\infty,T)\times\mathbb{R}^n \quad (38)$$

$$-\underbrace{\frac{1}{\lambda_\varphi \Phi'(1)}\sum_{i=1}^{n}c_i\int_0^{+\infty}\left\{\exp_q\left(-\lambda_\varphi\left(\Phi\left(\hat{\underline{\Psi}}_i\right)-1\right)\right)-1\right\}\nu(\mathrm{d}z_i)=0}_{\text{Jump effect considering risk and uncertainty}}$$

*with*

$$\hat{\underline{\Psi}}_i=\frac{\underline{\Psi}\left(t,\left\{x_j+\Delta_{i,j}z_j\right\}_{j=1,2,3,\ldots,n}\right)}{\underline{\Psi}\left(t,\left\{x_j\right\}_{j=1,2,3,\ldots,n}\right)} \tag{39}$$

*subject to the terminal condition*

$$\underline{\Psi}\left(T,\{x_i\}_{i=1,2,3,\ldots,n}\right)=\exp\left(p\sum_{i=1}^{n}x_i\right). \tag{40}$$

*In addition,* (38) *admits a smooth solution*

$$\underline{\Psi}=\exp\left(\sum_{i=1}^{n}\underline{\rho}_t^{(i)}x_i+\underline{\tau}_t\right) \tag{41}$$

*with time-dependent coefficients* $\left(\underline{\tau}_t,\left\{\underline{\rho}_t^{(i)}\right\}_{i=1,2,3,\ldots,n}\right)_{t\le T}$ *such that*

$$\underline{\rho}_t^{(i)}=\frac{1}{\dfrac{\sigma^2\left(1-\underline{\lambda}\right)}{2}+\left(\dfrac{1}{p}-\dfrac{\sigma^2\left(1-\underline{\lambda}\right)}{2}\right)e^{r_i(T-t)}} \tag{42}$$

*and*

$$\frac{\mathrm{d}\underline{\tau}_t}{\mathrm{d}t}+a\sum_{i=1}^{n}c_i\underline{\rho}_t^{(i)}-\frac{1}{\lambda_\varphi\Phi'(1)}\sum_{i=1}^{n}c_i\int_0^{+\infty}\left\{\exp_q\left(-\lambda_\varphi\left(\Phi\left(e^{\underline{\rho}_t^{(i)}z_i}\right)-1\right)\right)-1\right\}\nu(\mathrm{d}z_i)=0,\quad i=1,2,3,\ldots,n. \tag{43}$$

***Remark 3*** We have excluded the case $0<q<1$ from the assumption of **Proposition 2**, while the Kullback–Leibler case $q=1$ is allowed, in contrast to **Proposition 1**.

***Remark 4*** The exponential disutility of interest in applications is that of a stationary state, provided that the data considered do not exhibit significant trends. In this case, the first terms in the exponentials of (32) and (41) are omitted, and the second term is evaluated under the formal limit of $t\to-\infty$. This approach is described in **Section 4**.

***Remark 5*** One may be interested in the sample paths of supJCIR processes, but accurately simulating their sample paths is a difficult task at this stage. Indeed, it has been pointed out that even in the absence of jumps, common numerical methods for CIR processes converge arbitrarily slowly when $\sigma^2 r_i$ is significantly larger than $ac_i$ [64], which is typical in Markovian lifts [19]. However, the proposed approach based on the HJB equations can completely avoid this issue, as its resolution is free from statistical methods that rely on sample paths.

## 4. Case study

### 4.1 Study site

We apply the closed-form solutions obtained in **Propositions 1–2** to real data of environmental variables whose dynamics are identified as supJCIR processes, $Y$. The study site was taken to be the Kisuki Point, which is midstream of the main branch of the Hii River in Shimane Prefecture, Japan (**Figure 1**). The river reach containing the Kisuki point has been a key sampling point in previous research because it is considered to represent the nutrient loading point from the watershed of the Hii River [65,66]. The downstream reaches of the main branch of the Hii River include two brackish lakes, Lake Shinji and Lake Nakaumi, both of which are Lamar sites (https://www.ramsar.org/wetland/japan, accessed on June 17, 2023). Therefore, monitoring and assessing the water quality of the upstream reaches is of fundamental importance for sustainable conservation of the water environment and ecosystems in the Hii River.

We analyzed 30-year data of weekly sampled water quality indices at the Kisuki point until the end of 2021, with Total Nitrogen (TN, in the unit mg/L, from August 20, 1991 to December 28, 2021) and sulfate ions ($SO_4^{2-}$, in the unit mg/L, from March 2, 1993 to December 28, 2021 (**Figures 2(a)-(b)**). These are representative river environmental indicators of Hii River [67-69] possibly covering pollution sources; residential areas and agricultural activities such as fertilization for TN; fertilization and geological origin for $SO_4^{2-}$. The data had already been presented by Takeda [70]; however, no mathematical model had been fitted to them at Kisuki point. Therefore, a mathematical model was applied to this unique dataset for the first time as part of this study. We chose the two aforementioned indices owing to their contrasting properties—SDE of the former is sensitive to jumps, whereas that of the latter is insensitive to them. According to Takeda [70], no specific trend was detected in the water quality data at the Kisuki point; therefore, we applied the supJCIR process in a stationary state.

### 4.2 Parameter estimation

We assume that both diffusive and jump noises exist in the SDEs of the time series data and that the Lévy measure $v(\mathrm{dz})$ is that of an exponential distribution, as it yields a tractable jump model:

$$v(\mathrm{dz}) = \underbrace{\mu}_{\text{Jump intensity}} \times \underbrace{\beta \exp(-\beta z)\mathrm{dz}}_{\text{Jump size probability density}}, \quad z > 0 \tag{44}$$

with parameters $\beta, \mu > 0$. Now, the parameters to be identified are $\theta, \omega, \mu, \beta, \sigma, a$. Elementary calculation yields stationary moments:

$$\text{Average} = R\left(a + M_1\right) = R\left(a + \frac{\mu}{\beta}\right), \tag{45}$$

$$\text{Variance} = R\left(\frac{M_2}{2} + \frac{\sigma^2}{2}\left(a + M_1\right)\right) = R\left(\frac{\mu}{\beta^2} + \frac{\sigma^2}{2}\left(a + \frac{\mu}{\beta}\right)\right), \tag{46}$$

$$\text{Skewness} = \frac{R\left(\frac{\sigma^4}{2}\left(a + M_1\right) + \frac{\sigma^2 M_2}{2} + \frac{M_3}{3}\right)}{\text{Vartiance}^{3/2}} = \frac{R\left(\frac{\sigma^4}{2}\left(a + \frac{\mu}{\beta}\right) + \frac{\sigma^2\mu}{\beta^2} + \frac{2\mu}{\beta^3}\right)}{\text{Vartiance}^{3/2}}. \tag{47}$$

Model parameters of each water quality index are identified using a two-step method originally proposed in [50] for purely jump-driven processes. We modify the second step to adapt it to our model. In the first step, the parameters $\theta, \omega$ of the probability measure $\pi$ are fitted using a naïve least-squares method on the empirical and theoretical autocorrelation functions. In the second step, based on the identified values of $\theta, \omega$, the remaining parameters $\mu, \beta, \sigma, a$ are identified so that the sums of relative errors of the average, variance, and skewness are minimized. Optimization is performed to minimize the following error metric, where the subscripts "e" and "m" represent empirical and modelled quantities, respectively:

$$\left(\frac{\text{Average}_\text{e} - \text{Average}_\text{m}}{\text{Average}_\text{e}}\right)^2 + \left(\frac{\text{Variance}_\text{e} - \text{Variance}_\text{m}}{\text{Variance}_\text{e}}\right)^2 + \left(\frac{\text{Skewness}_\text{e} - \text{Skewness}_\text{m}}{\text{Skewness}_\text{e}}\right)^2. \tag{48}$$

We add a constraint in the second step to mitigate the risk of the minimization problem becoming underdetermined because of the four unknown parameters, whereas the error metric (48) is the sum of the three normalized statistics. We assume that $100y$ % ($0 \le y \le 1$) of the average is explained by the diffusive part; hence, the remaining $100(1-y)$ % is explained by the jump part. By (45), this implies that:

$$(1-y)a = y\frac{\mu}{\beta} \text{ and hence } \mu = \frac{1-y}{y}\beta a. \tag{49}$$

We then prescribe $y$ and enforce (49) while minimizing (48). The goodness-of-fit of the identified model with respect to the empirical data is examined by comparing not only the statistics but also the probability density functions.

**Tables 1** presents the estimated parameter values of $\theta, \omega$ for TN and $SO_4^{2-}$. **Tables 2–3** present the estimated parameter values of $\mu, \beta, \sigma, a$ for TN and $SO_4^{2-}$ at different levels of $y$. **Tables 4–5** present the relative errors between each pair of empirical and theoretical moments. **Table 1** suggests that both the time series data of TN and $SO_4^{2-}$ exhibit long memories such that the autocorrelation function is not integrable in $(0, +\infty)$ [27,50]. Hence, it is essentially subexponential. This implies that the use of the supJCIR process is critical in our application. Further, **Figures 3(a)–(b)** depict the fitted autocorrelation functions for both indices, demonstrating the reasonable fit of the theoretical model.

**Tables 2–3** indicate that the model of TN is more sensitive to jumps than that of $SO_4^{2-}$, suggesting that the time series of $SO_4^{2-}$ is almost exclusively driven by diffusive dynamics rather than jumps. Indeed, the jump intensity $\mu$ has been identified to be significantly larger than 1 and the mean jump size $\beta^{-1}$ significantly smaller than 1, under which jumps are possibly approximated by continuous Brownian noises [71,72]. **Tables 4–5** demonstrate that the empirical moments are comparably accurate among the theoretical models for both indices. In contrast, as illustrated in **Figures 4(a)–(b)**, accounting for jumps affects the probability density functions of TN. **Figure 4(a)** demonstrates that the value around $y = 0.95$ yields a reasonable fit between the empirical and theoretical probability density functions. In **Figure 4(b)**, the average of the empirical $SO_4^{2-}$ is underestimated by the theoretical one. Moreover, for $SO_4^{2-}$, we also examine model identification without considering the last term of (48), which is also plotted in **Figure 4(b)** (see also **Table 5**). The theoretical model constructed without considering the last term of (48) exhibits a good fit with the empirical one, implying that considering a larger number of statistics does not always improve model performance. Similar reduced error metric is used to identify the model of TN—no significant differences from the data presented in **Table 4** is observed. Thus, the results of this application are not reported here.

**Table 1.** Estimated $\theta, \omega$ for TN and $SO_4^{2-}$.

| | TN | $SO_4^{2-}$ |
|---|---|---|
| $\theta$ (1/day) | 0.715 | 0.163 |
| $\omega$ (-) | 1.62 | 1.52 |

**Table 2.** Estimated $\mu, \beta, \sigma, a$ for TN at different levels of $y$.

| $y$ | 1 | 0.999 | 0.99 | 0.95 | 0.94 |
|---|---|---|---|---|---|
| $\mu$ (1/day) | 0 | 0.0000188 | 0.000587 | 0.00644 | 0.00843 |
| $\beta$ (L/mg) | - | 0.0706 | 0.221 | 0.484 | 0.528 |
| $\sigma$ (mg$^{1/2}$/L$^{1/2}$) | 0.574 | 0.504 | 0.438 | 0.276 | 0.236 |
| $a$ (1/day) | 0.238 | 0.266 | 0.263 | 0.253 | 0.250 |

**Table 3.** Estimated $\mu, \beta, \sigma, a$ for $SO_4^{2-}$ at different levels of $y$. Results for the model constructed without considering the last term of (48) with $y = 1$ ("Reduced" in the table) are also reported.

| $y$ | 1 | 0.999 | 0.99 | 0.95 | 0.94 | Reduced |
|---|---|---|---|---|---|---|
| $\mu$ (1/day) | 0 | 1.52E+06 | 1.11E+03 | 7.59.E+07 | 9.11.E+07 | 0 |
| $\beta$ (L/mg) | - | 5.09E+09 | 3.73E+05 | 5.09.E+09 | 5.09.E+09 | - |
| $\sigma$ (mg$^{1/2}$/L$^{1/2}$) | 1.75 | 1.75 | 1.75 | 1.75 | 1.75 | 1.44 |
| $a$ (1/day) | 0.298 | 0.298 | 0.295 | 0.284 | 0.281 | 0.403 |

**Table 4.** Comparison between each pair of empirical and theoretical moments for TN at different levels of $y$.

| | | Model for different $y$ | | | | |
|---|---|---|---|---|---|---|
| | Empirical | 1 | 0.999 | 0.99 | 0.95 | 0.94 |
| Average | 6.01E-01 | 5.38E-01 | 6.01E-01 | 6.01E-01 | 6.01E-01 | 6.01E-01 |
| Variance | 8.49E-02 | 8.87E-02 | 8.49E-02 | 8.49E-02 | 8.49E-02 | 8.49E-02 |
| Skewness | 1.06E+01 | 1.11E+00 | 1.06E+01 | 1.06E+01 | 1.06E+01 | 1.06E+01 |

**Table 5.** Comparison between each pair of empirical and theoretical moments for TN at different levels of $y$. Results for the model constructed without considering the last term of (48) with $y = 1$ ("Reduced" in the table) are also reported.

| | | Model for different $y$ | | | | | |
|---|---|---|---|---|---|---|---|
| | Empirical | 1 | 0.999 | 0.99 | 0.95 | 0.94 | Reduced |
| Average | 4.76E+00 | 3.53E+00 | 3.53E+00 | 3.53E+00 | 3.53E+00 | 3.53E+00 | 4.76E+00 |
| Variance | 4.96E+00 | 5.40E+00 | 5.40E+00 | 5.40E+00 | 5.40E+00 | 5.40E+00 | 4.96E+00 |
| Skewness | 1.78E+00 | 1.32E+00 | 1.32E+00 | 1.32E+00 | 1.32E+00 | 1.32E+00 | 9.36E-01 |

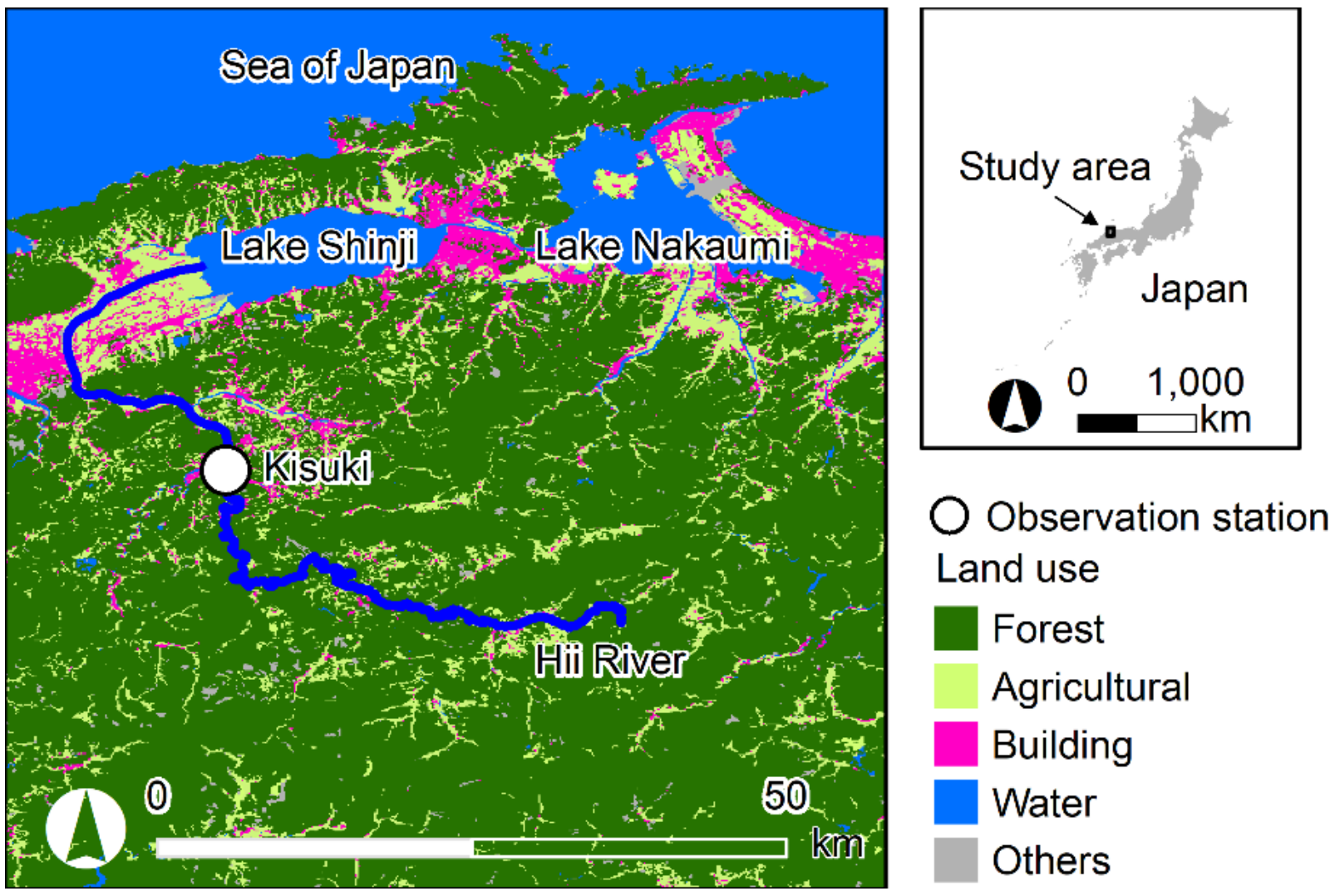

**Figure 1.** Map of the study site.

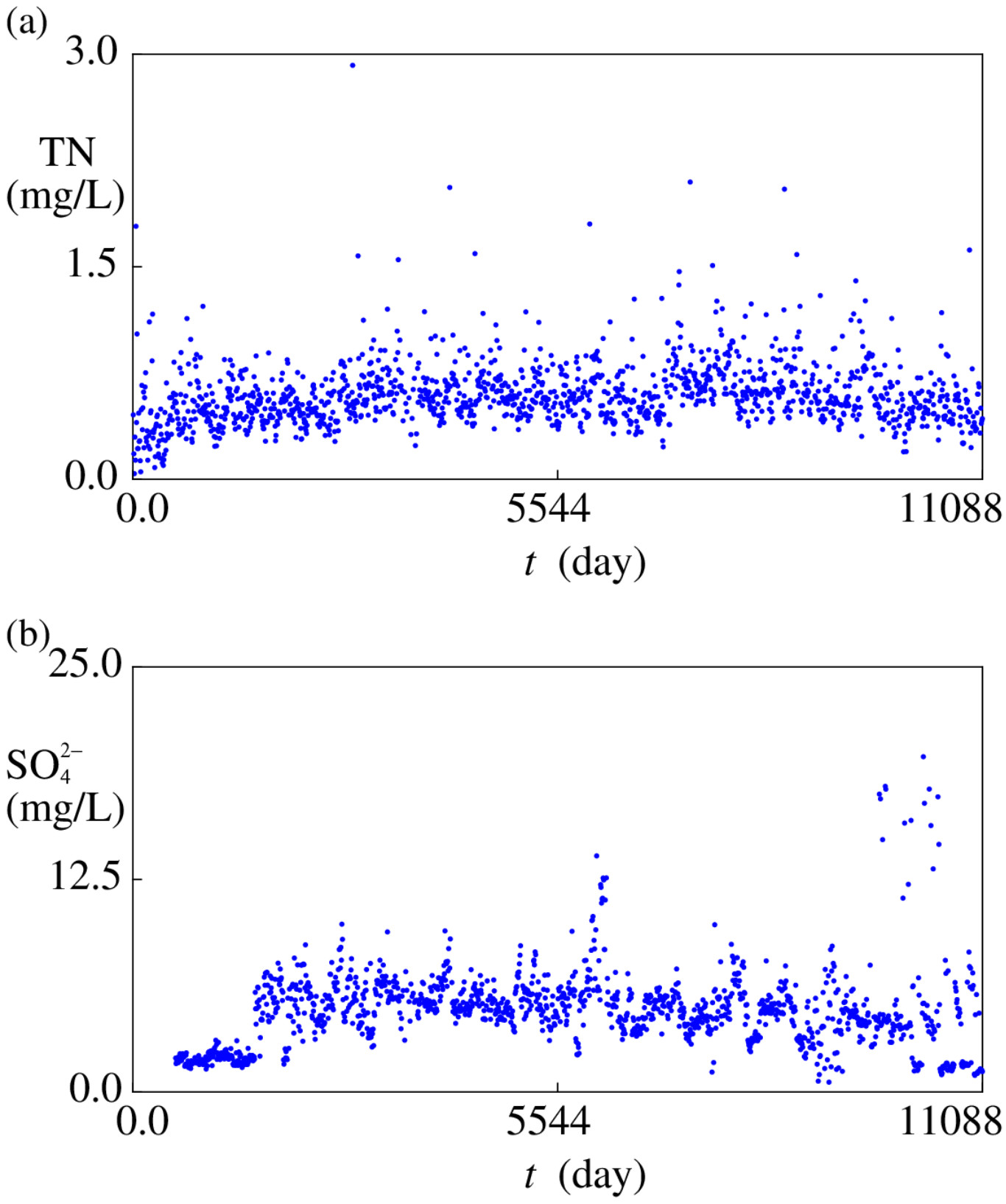

**Figure 2.** Time series data of (a) TN and (b) $SO_4^{2-}$. The 0th (day) is August 20, 1991.

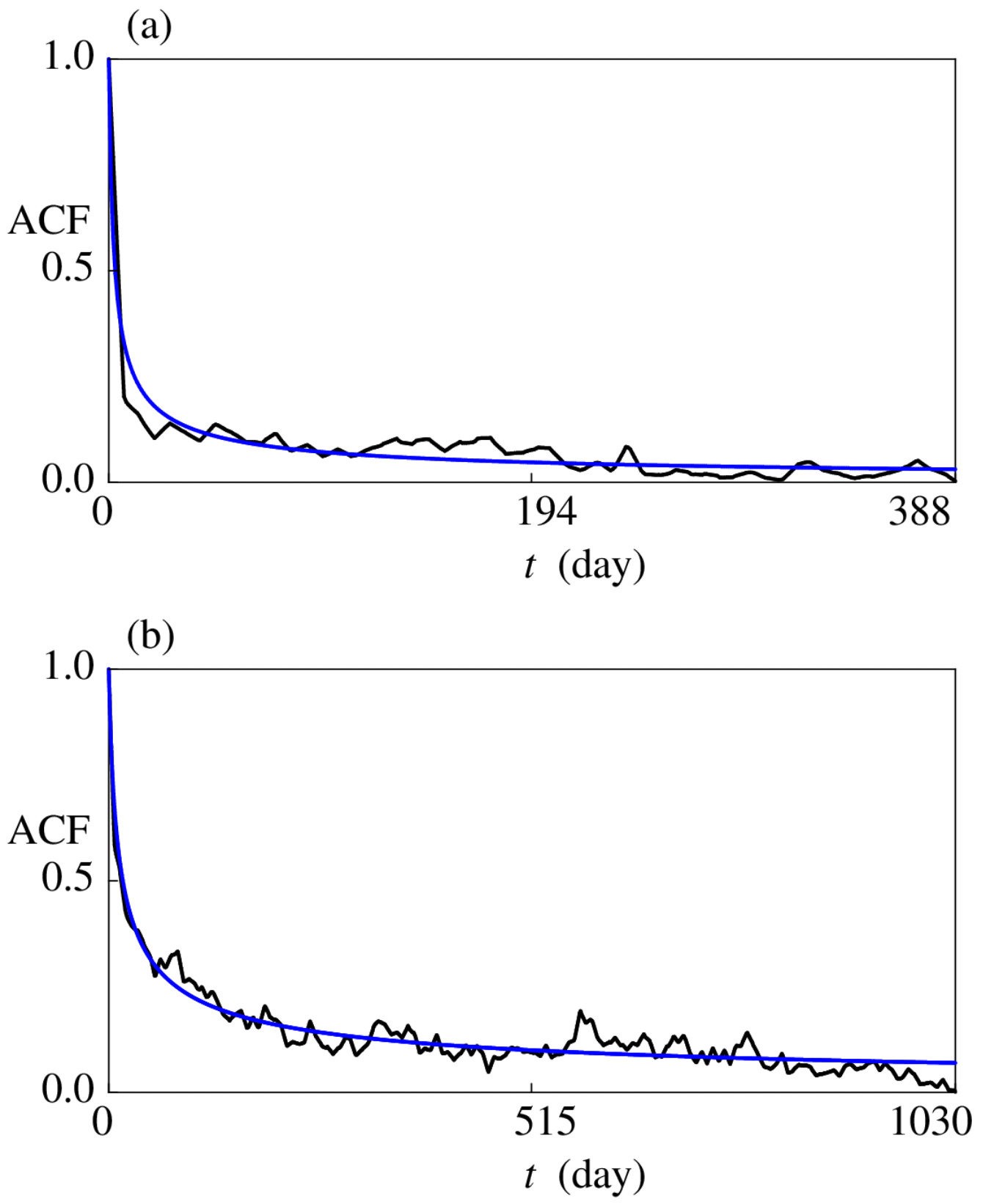

**Figure 3.** Empirical (black) and theoretical (blue) autocorrelation functions (ACFs) of (a) TN and (b) $SO_4^{2-}$.

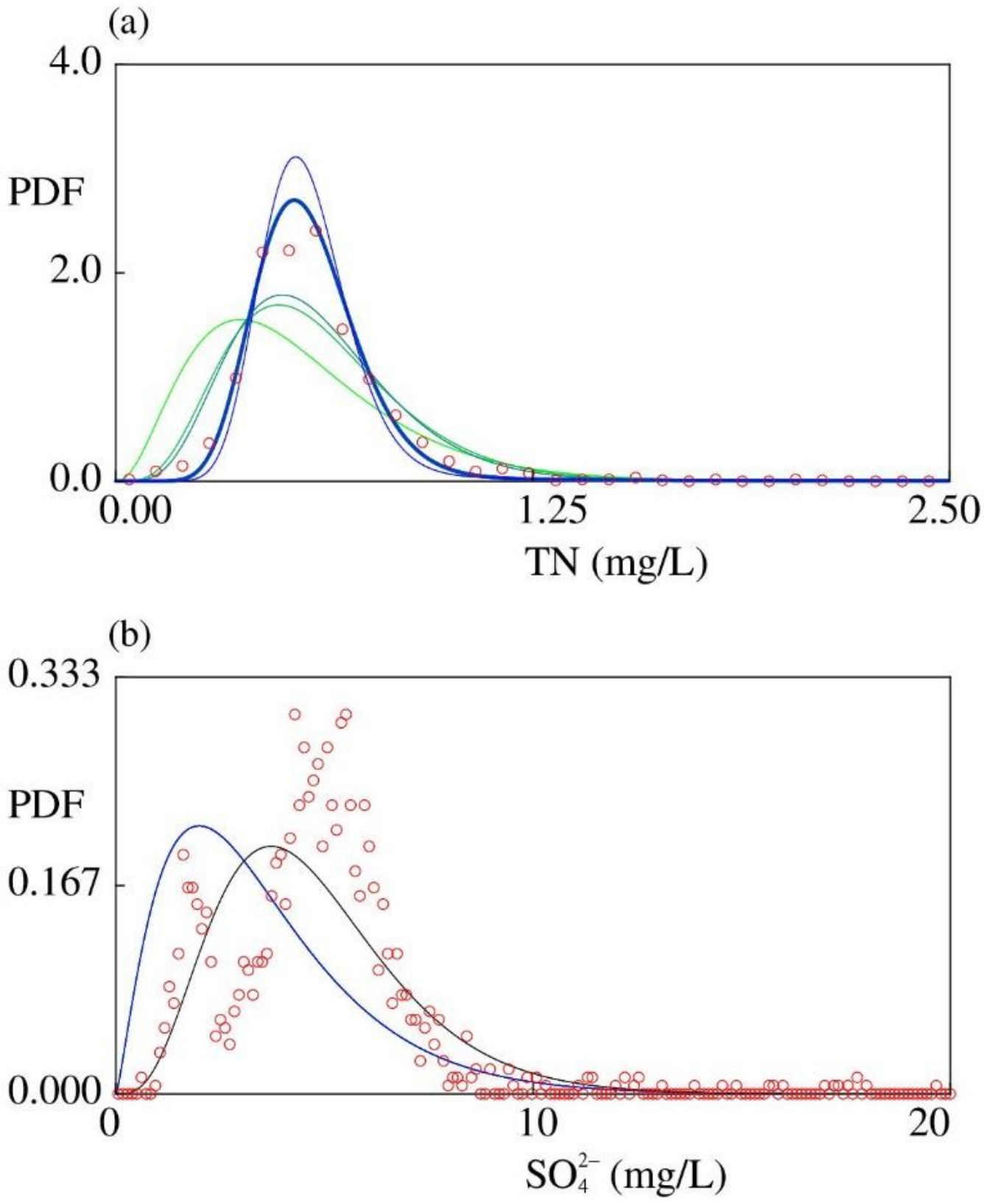

**Figure 4.** Empirical (circles) and theoretical (curves) probability density functions (PDFs) of (a) TN and **(b)** $SO_4^{2-}$. In the panel (a), the model of the theoretical PDFs moves from the left to the right as $y$ decreases along 1, 0.999, 0.99, 0.95, and 0.94. The theoretical PDF with $y = 0.95$ is plotted using a bold curve. In panel (b), all theoretical PDFs identified based on (48) are almost the same and are overlapping (See also **Table 5**), which are highlighted in blue. In panel (b), the black curve corresponds to the model identified not considering the last term of (48).

### 4.3 Application

For TN, we compute the exponential disutilities using the supJCIR process, assuming a stationary state based on the closed-form solutions given in **Propositions 1–2** combined with **Remark C1** in **Appendix C**. Qualitatively, the same computational results are obtained for $SO_4^{2-}$, which is observed to follow a simpler model without jumps. The results of $SO_4^{2-}$ are presented in **Appendix D**.

We examine different parameter values of $p, \lambda_\phi, \lambda_\varphi$ and the normalized exponential utilities such that they are equal to 1 when there is no misspecification.

$$U = \frac{\bar{\Psi}\big|_{\text{With uncertainty}}}{\bar{\Psi}\big|_{\text{No uncertainty}}} \in [1, +\infty) \text{ (for upper bound) or } U = \frac{\underline{\Psi}\big|_{\text{With uncertainty}}}{\underline{\Psi}\big|_{\text{No uncertainty}}} \in (0, 1] \text{ (for lower bound).} \quad (50)$$

The disutility with no uncertainty can be obtained by $\lambda_\phi, \lambda_\varphi \to +0$. The benchmark model for TN is selected as the one with $y = 0.95$.

**Figures 5–6** depict the normalized upper and lower bounds of the exponential distributions for TN with $\Phi = x$ and different values of $p$, respectively. Here, we consider $q = 0.75$ as the upper bound and $q = 1.25$ as the lower bound; hence, only the influences of the model uncertainties are evaluated. We examine the parameter values $(\lambda_\phi, \lambda_\varphi) = (0.01 + 2000k, 0.01 + 200l)$ for the upper bound and $(\lambda_\phi, \lambda_\varphi) = (0.01 + 5000k, 0.01 + 200l)$ for the lower bound ($k, l = 0, 1, 2, ..., 20$). Comparing **Figures 5(a)–(b)** for the upper bound reveals that the dependence of the disutilities on the weights $(\lambda_\phi, \lambda_\varphi)$ is qualitatively identical to that on the different values of $p$, whereas the disutilities corresponding to larger values of $p$ are more sensitive to model uncertainty. The same applies to the case of lower bounds. Therefore, exponential disutility should be designed by considering the risk of a larger accidental increase in TN in the watershed of the Hii River due to different factors, e.g., industrial pollution from sewage treatment plants or paddy fields with intensive fertilization during spring.

**Figures 7–8** present the normalized upper bound with $\Phi = x^{3/2}$ and lower bound with $\Phi = x^{1/2}$ of the exponential distributions for TN with respect to different values of $p$, respectively. The Orlicz risk uses nonlinear $\Phi$. The influence of the nonlinear and convex $\Phi = x^{3/2}$ for the upper bound is clearly indicated in **Figure 7**, with a sharper increase, and hence larger model uncertainty, in disutility for larger $(\lambda_\phi, \lambda_\varphi)$ than for $\Phi = x$ in **Figure 5**. The use of the Orlicz risk in this case yields a more sensitive risk indicator for severe environmental pollution, which may not necessarily be suitable for evaluating environmental risks caused by a moderate increase in TN. In the study area, such a case may arise in the future if urbanization continues around the river or due to land use in the watershed transitions to a larger area of grazing land. In contrast, the lower bound case depicted in **Figure 8** indicates that the use of the concave $\Phi = x^{1/2}$ does not affect the exponential disutilities significantly, suggesting that the use of the affine $\Phi = x$ suffices for risk evaluation in this case.

We also analyze impact of the parameter $q$ on the modeling of relative entropy concerning jumps in TN. **Figures 9–10** depict the normalized upper bound with $q = 0.5$ and lower bound with $q = 1.5$ of the exponential disutilities for TN with respect to different values of $p$, respectively. The parameter $q$ is selected to be more distant from the Kullback–Leibler case $q = 1$ than that in **Figures 7–8**. **Figure 9** demonstrates that the very high sensitivity of the upper-bound of the exponential disutility is mitigated by choosing smaller $q$, suggesting that the choice of $q$ and hence the relative entropy in Orlicz risk is important during the design of the upper bound. In contrast, **Figure 10** reveals that the use of a larger $q$ does not affect the normalized lower bound of the disutility, demonstrating its robustness with respect to relative entropy.

The obtained upper and lower bounds of the normalized exponential disutilities of TN suggest that Orlicz risks can be employed as a flexible mathematical tool for the evaluation of environmental pollution. Based on model application, the lower bound is observed to be insensitive to its functional form, at least for the parameter values considered in this study. As TN is an indicator of environmental pollution including eutrophication, an optimistic lower bound may be less important than a pessimistic upper bound. Nevertheless, some water quality indices, such as dissolved oxygen, should be maintained above a certain threshold to preserve habitat suitability for a variety of aquatic species, such as riverine fish species [73,74]. The lower bound is suitable in such cases.

Finally, we discuss the influence of model uncertainties on the moments and the autocorrelation function. In the stationary state, according to **Propositions 1–2** and **Remark B1** in **Appendix B**, the SDE (9) is distorted to

$$\mathrm{d}X_t^{(i)} = \left(ac_i - r_i\left(1-\xi\right)X_t^{(i)}\right)\mathrm{d}t + \sigma\sqrt{r_i X_t^{(i)}}\mathrm{d}W_t^{(i)} + \mathrm{d}\bar{J}_t^{(i)} \tag{51}$$

if the exponential disutility is overestimated using Orlicz risk, and to

$$\mathrm{d}X_t^{(i)} = \left(ac_i - r_i\left(1+\xi\right)X_t^{(i)}\right)\mathrm{d}t + \sigma\sqrt{r_i X_t^{(i)}}\mathrm{d}W_t^{(i)} + \mathrm{d}\underline{J}_t^{(i)} \tag{52}$$

if it is underestimated. Here, $\xi = \lambda_\varphi \Phi'(1)\sigma^2 p > 0$ is a constant, and $\bar{J}$ (resp., $\underline{J}$ ) denotes a subordinators with the Lévy measure $\exp_q\left\{\lambda_\varphi \Phi\left(e^{pz}\right) - \lambda_\varphi\right\}v(\mathrm{d}z)$ (resp., $\exp_q\left\{-\lambda_\varphi\left(\Phi\left(e^{pz}\right)-1\right)\right\}v(\mathrm{d}z)$ ). Both (51) and (52) are JCIR processes, and hence their moments and autocorrelation functions can be explicitly found. We focus on the autocorrelation function because it is directly related to the long memory. For statistical moments, see **Appendix D**.

**Figure 11** compares the autocorrelation functions in the case without any uncertainty in the upper-bound case (red), and the lower-bound case (green) (For their formulae, see **Remark B1)**. In the upper-bound case (resp., lower-bound case), model uncertainty leads to longer memory (resp., shorter memory). In the upper-bound case, longer memory implies the expectation of prolonged accidental environmental pollution events. From this perspective, environmental risk assessment based on Orlicz risk for the upper bound leads to more pessimistic prediction of both the size and duration of pollution. It should also be noted that according to **Remark B1**, the power of the decaying speed of the autocorrelation does not change from that in (18) under the worst-case models, demonstrating a structural robustness of the supJCIR process evaluated through the Orlicz risks.

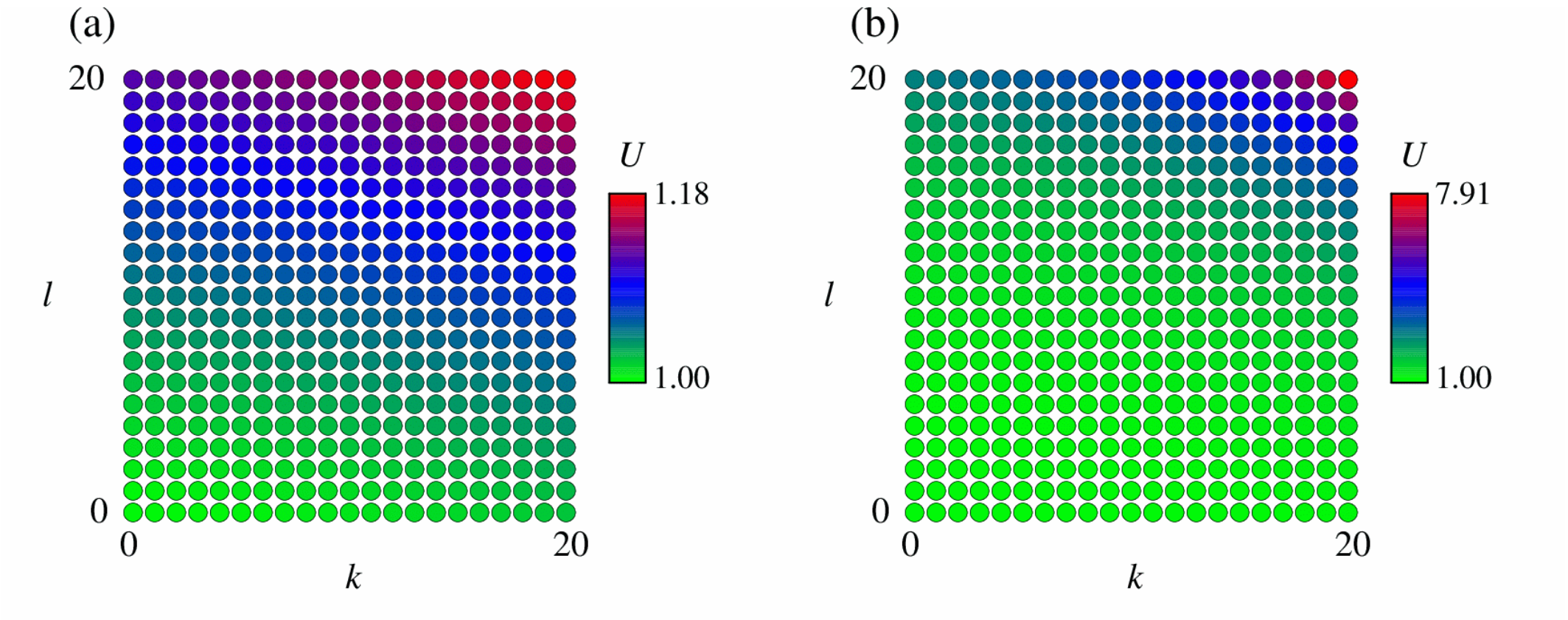

**Figure 5.** The upper bound of the normalized exponential disutility $U$ with (a) $p=0.02$ and (b) $p=0.04$, where $\Phi=x$, $q=0.75$, and $\left(\lambda_{\phi},\lambda_{\varphi}\right)=\left(0.01+2000k,0.01+200l\right)$ ($k,l=0,1,2,...,20$). The exponential disutilities in the case without model uncertainty are (a) 1.11 and (b) 1.24.

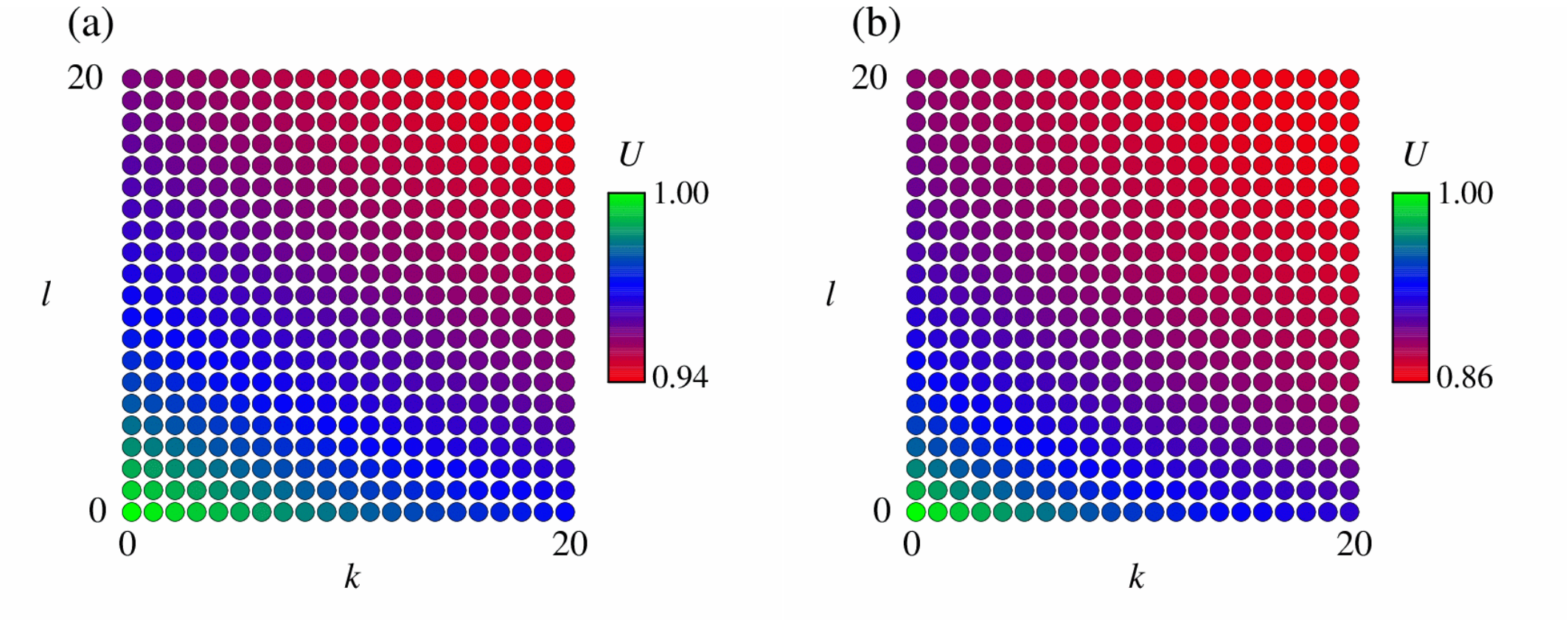

**Figure 6.** The lower bound of the normalized exponential disutility $U$ with (a) $p=0.02$ and (b) $p=0.04$, where $\Phi=x$, $q=1.25$, and $\left(\lambda_{\phi},\lambda_{\varphi}\right)=\left(0.01+5000k,0.01+200l\right)$ ($k,l=0,1,2,...,20$). The exponential disutilities in the case without model uncertainty are (a) 1.11 and (b) 1.24.

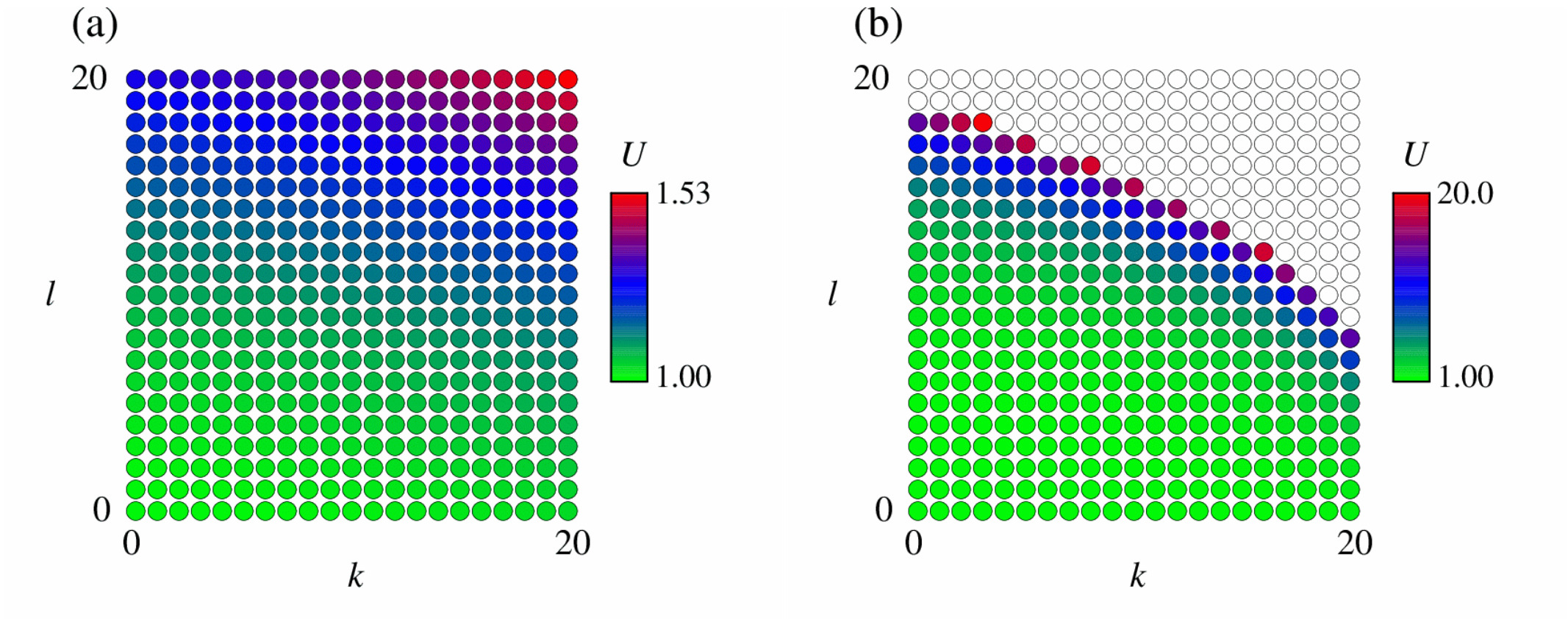

**Figure 7.** The upper bound of the normalized exponential disutility $U$ with (a) $p = 0.02$ and (b) $p = 0.04$, where $\Phi = x^{3/2}$ and the other conditions are identical to those in **Figure 5**. The exponential disutilities in the case without model uncertainty are (a) 1.11 and (b) 1.25. The white plots in the panel (b) represent the cases where the exponential disutility diverges.

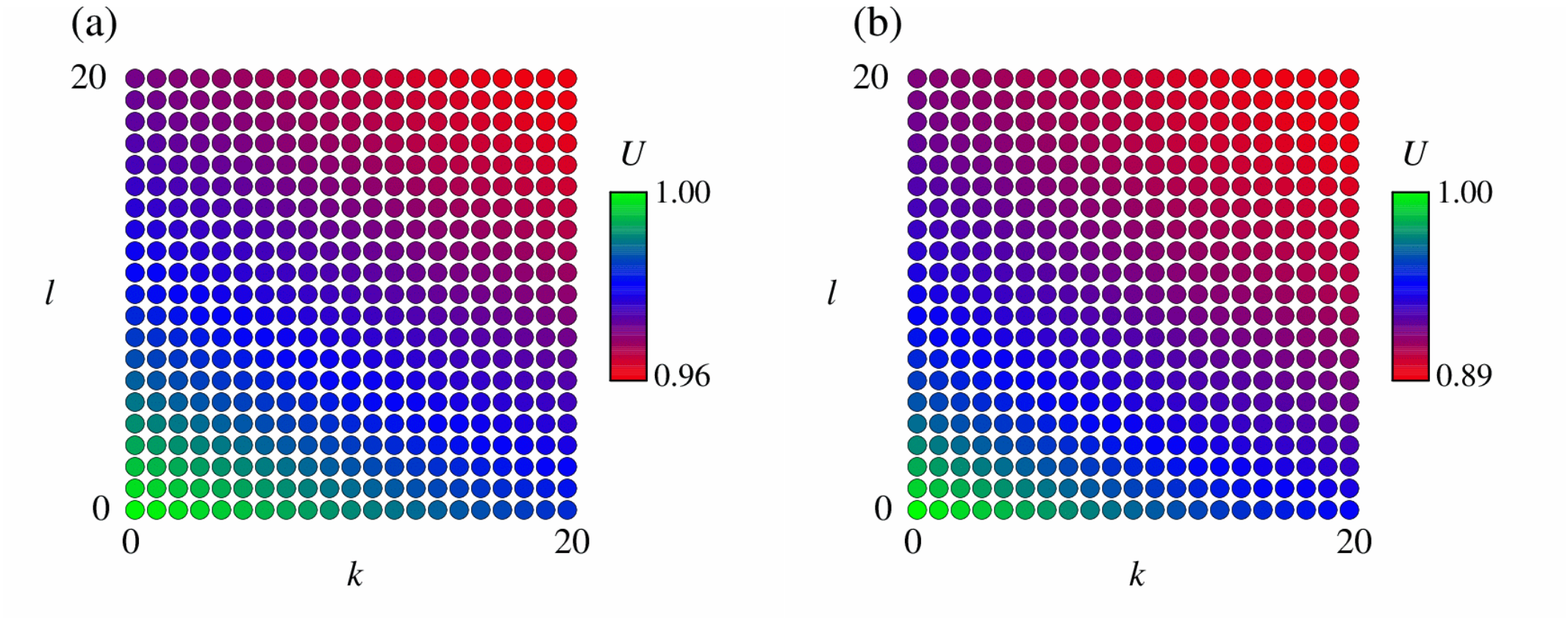

**Figure 8.** The lower bound of the normalized exponential disutility $U$ with (a) $p = 0.02$ and (b) $p = 0.04$, where $\Phi = x^{1/2}$ and the other conditions are identical to those in **Figure 6**. The exponential disutilities in the case without model uncertainty are (a) 1.11 and (b) 1.24.

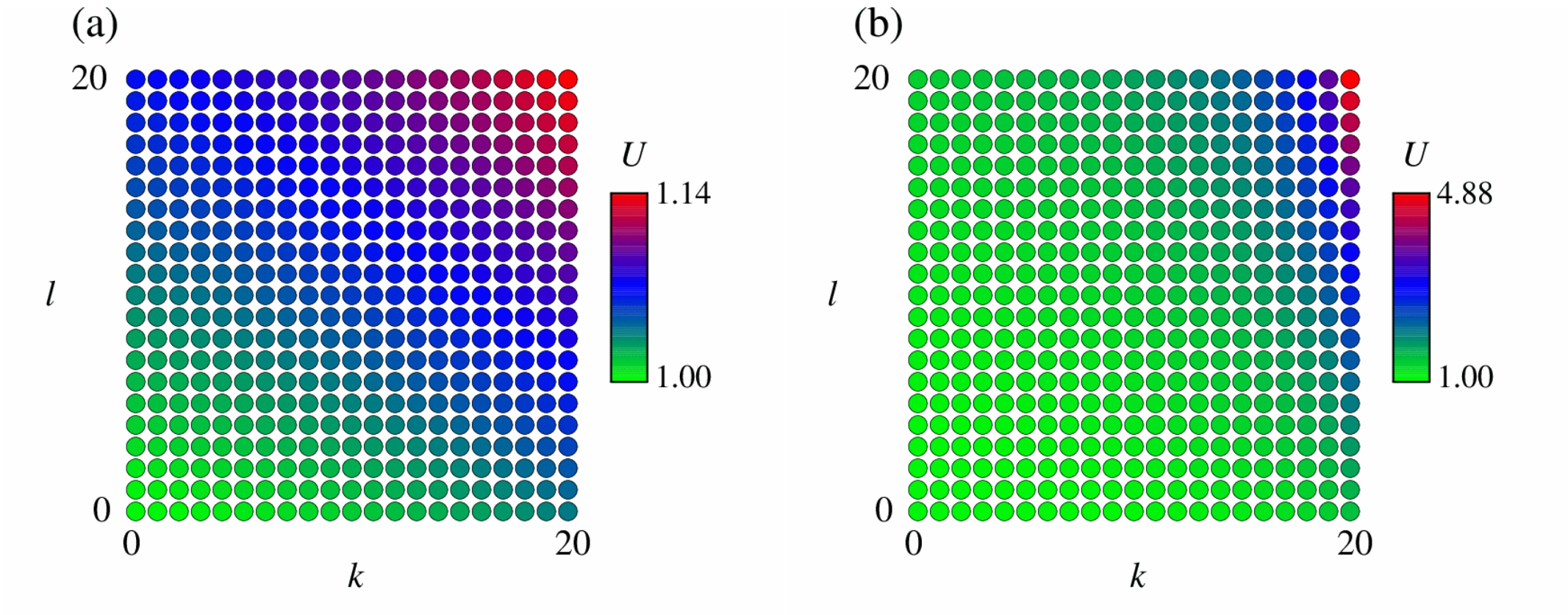

**Figure 9.** The upper bound of the normalized exponential disutility $U$ with (a) $p = 0.02$ and (b) $p = 0.04$, where $q = 0.5$ and the other conditions are the same as in **Figure 7**. The exponential disutilities in the case without model uncertainty are (a) 1.11 and (b) 1.25.

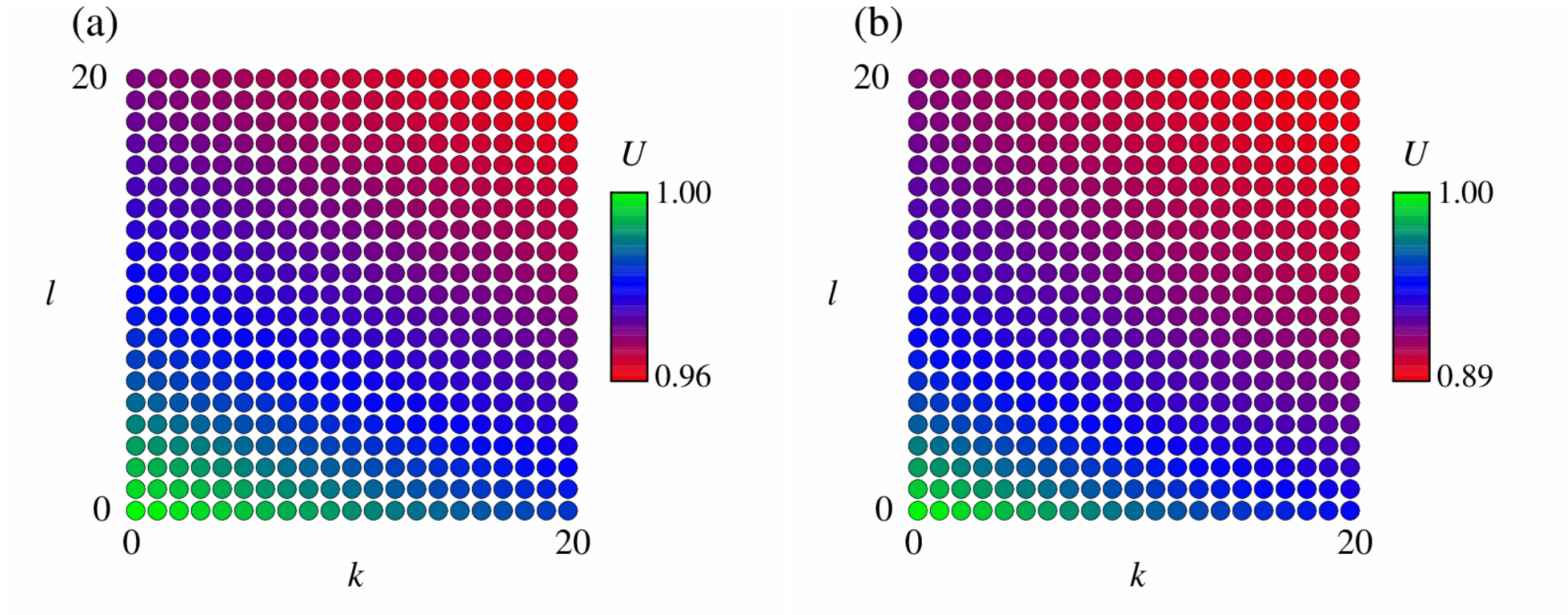

**Figure 10.** The lower bound of the normalized exponential disutility $U$ with (a) $p = 0.02$ and (b) $p = 0.04$, where $q = 1.5$ and the other conditions are the same as in **Figure 8**. The exponential disutilities in the case without model uncertainty are (a) 1.11 and (b) 1.24.

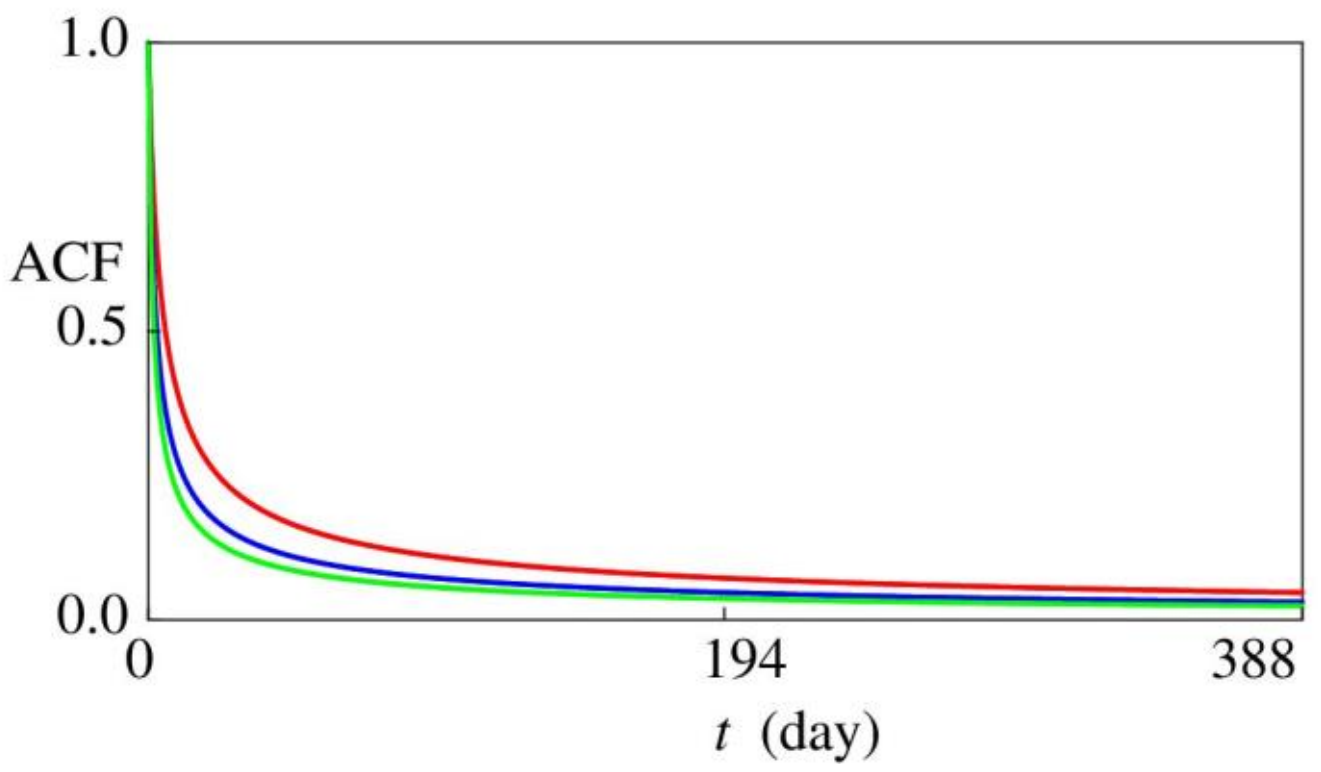

**Figure 11.** Autocorrelation functions (ACFs) in cases with no uncertainty (blue), upper bound (red), and lower bound (green) of TN.

## 5. Conclusions

We proposed a recursive (dis)utility to evaluate the risk of extreme values of environmental variables as long-memory processes. The core of our methodology is the Orlicz risks, from which recursive utilities are explicitly computed using HJB equations. The well-defined HJB equations were analyzed and applied to real-world water quality data of a Japanese river.

The use of exponential disutility is a novel approach, as it has not yet been combined with Orlicz risks, unlike other major types of disutilities, such as power and logarithmic ones. The latter cases were not considered in this study because no analytical solutions to the corresponding HJB equations were found. A purely numerical approach that can deal with high-dimensional HJB equations [75] is necessary to resolve this technical issue. Another limitation of our mathematical framework was the specific form of the SDE. Theoretically, both the diffusion and jump parts can be made more sophisticated so that the resulting SDE better represents the target time series data.

This study demonstrated that Orlicz risks can be applied not only to the SDEs considered in this study but also to more general jump-diffusion SDEs if the corresponding HJB equation can be solved analytically or efficiently via numerical computation. Markovian lifts often applied to Volterra-type SDEs; hence, the finite-dimensional approach employed in this study is sufficiently general and can be adapted to their risk analysis. The proposed mathematical framework can also be applied to affine self-exciting jump processes that describe discharge time-series data [76]. The authors are currently considering the application of Orlicz risk to green renewable power generation. Its resolution has been reduced to solving an HJB equation where the underlying SDE is neither linear nor affine. Coupling with the hydrodynamic will also be interesting [77].

**Acknowledgments** The authors would like to express their gratitude towards Dr. Ikuo Takeda of Shimane University for providing time-series data on the water quality indices of Kisuki point, Hii River, Japan.

**Appendices of "Assessing fluctuations of long-memory environmental variables based on the robustified dynamic Orlicz risk"**

## Appendix A. Proof of Proposition 1

This proof proceeds via two steps. In the first step, the HJB equation is derived from the Orlicz risk. In the second step, a closed-form solution is obtained by directly guessing the solution.

Recall that $\Phi(1)=1$. The Orlicz risk for any $t<T$ with a sufficiently small $k>0$ can be rewritten as

$$\sup_{\mathbb{Q}(\phi,\varphi)}\left\{\begin{array}{l}\mathbb{E}_{\mathbb{Q}(\phi,\varphi)}\left[\Phi\left(\dfrac{\bar{\Psi}_{t+k}}{\bar{\Psi}_t}\right)\middle|\mathcal{F}_t\right]-\Phi(1)\\ -\left(\begin{array}{l}\dfrac{1}{\lambda_\phi}\mathbb{E}_{\mathbb{Q}(\phi,\varphi)}\left[\dfrac{1}{2}\displaystyle\sum_{i=1}^{n}\int_t^{t+k}\left(\phi_s^{(i)}\right)^2\mathrm{d}s\right]\\ +\dfrac{1}{\lambda_\varphi}\mathbb{E}_{\mathbb{Q}(\phi,\varphi)}\left[\displaystyle\sum_{i=1}^{n}\int_t^{t+k}\int_0^{+\infty}\dfrac{1}{1-q}\left(1-q-\left(\varphi_s^{(i)}(z)\right)^q+q\varphi_s^{(i)}(z)\right)v(\mathrm{d}z)\mathrm{d}s\right]\end{array}\right)\end{array}\right\}=0. \tag{53}$$

Assume $\{x_i\}_{i=1,2,3,\ldots,n}=\left\{X_t^{(i)}\right\}_{i=1,2,3,\ldots,n}$. We divide both sides by $k$, and later let $k\to+0$. Using the classical mean value theorem, we obtain (we omit the subscript representing time when there will be no confusion)

$$\begin{aligned}&\lim_{k\to+0}\frac{1}{k}\left\{\begin{array}{l}\dfrac{1}{\lambda_\phi}\mathbb{E}_{\mathbb{Q}(\phi,\varphi)}\left[\dfrac{1}{2}\displaystyle\sum_{i=1}^{n}\int_t^{t+k}\left(\phi_s^{(i)}\right)^2\mathrm{d}s\right]\\ +\dfrac{1}{\lambda_\varphi}\mathbb{E}_{\mathbb{Q}(\phi,\varphi)}\left[\displaystyle\sum_{i=1}^{n}\int_t^{t+k}\int_0^{+\infty}\dfrac{1}{1-q}\left(1-q-\left(\varphi_s^{(i)}(z)\right)^q+q\varphi_s^{(i)}(z)\right)v(\mathrm{d}z)\mathrm{d}s\right]\end{array}\right\}\\ &=\frac{1}{\lambda_\phi}\frac{1}{2}\sum_{i=1}^{n}\left(\phi^{(i)}\right)^2+\frac{1}{\lambda_\varphi}\sum_{i=1}^{n}\int_0^{+\infty}\frac{1}{1-q}\left(1-q-\left(\varphi^{(i)}(z)\right)^q+q\varphi^{(i)}(z)\right)v(\mathrm{d}z)\end{aligned}. \tag{54}$$

We also have

$$\Phi\left(\frac{\bar{\Psi}_{t+k}}{\bar{\Psi}_t}\right)=\Phi\left(\frac{\bar{\Psi}\left(t+k,\left\{X_{t+k}^{(i)}\right\}_{i=1,2,3,\ldots,n}\right)}{\bar{\Psi}\left(t,\left\{X_t^{(i)}\right\}_{i=1,2,3,\ldots,n}\right)}\right)=\Phi\left(\frac{\bar{\Psi}\left(t+k,\left\{X_{t+k}^{(i)}\right\}_{i=1,2,3,\ldots,n}\right)}{\bar{\Psi}\left(t,\{x_i\}_{i=1,2,3,\ldots,n}\right)}\right). \tag{55}$$

By applying Itô's formula to $\Xi\left(t+k,\left\{X_{t+k}^{(i)}\right\}_{i=1,2,3,\ldots,n}\right)=\Phi\left(\dfrac{\bar{\Psi}\left(t+k,\left\{X_{t+k}^{(i)}\right\}_{i=1,2,3,\ldots,n}\right)}{\bar{\Psi}\left(t,\left\{X_t^{(i)}\right\}_{i=1,2,3,\ldots,n}\right)}\right)$ as a function of $\left(t+k,\left\{X_{t+k}^{(i)}\right\}_{i=1,2,3,\ldots,n}\right)$, we have

$$
\begin{aligned}
&\Xi\left(t+k,\left\{X_{t+k}^{(i)}\right\}_{i=1,2,3,\ldots,n}\right)-\Xi\left(t,\left\{X_{t}^{(i)}\right\}_{i=1,2,3,\ldots,n}\right) \\
&=\int_t^{t+k}\left(\begin{array}{l}\dfrac{\partial \Xi_s}{\partial s}+\displaystyle\sum_{i=1}^{n}\left(\left(ac_i-r_iX_s^{(i)}+\sigma\sqrt{r_iX_s^{(i)}}\phi_s^{(i)}\right)\dfrac{\partial \Xi_s}{\partial X_s^{(i)}}+\dfrac{\sigma^2}{2}r_iX_s^{(i)}\dfrac{\partial^2 \Xi_s}{\partial X_s^{(i)}\partial X_s^{(i)}}\right) \\ +\displaystyle\sum_{i=1}^{n}c_i\int_0^{+\infty}\varphi_s^{(i)}\left(z_i\right)\left\{\Xi\left(s,\left\{X_s^{(j)}+\Delta_{i,j}z_j\right\}_{j=1,2,3,\ldots,n}\right)-\Xi_s\right\}v\left(\mathrm{d}z_i\right)\end{array}\right)\mathrm{d}s \\
&\quad+\int_t^{t+k}\sum_{i=1}^{n}\sigma\sqrt{r_iX_s^{(i)}}\frac{\partial \Xi_s}{\partial X_s^{(i)}}\mathrm{d}W_s^{(i)}+\int_t^{t+k}\sum_{i=1}^{n}\sum_{s\in(t,t+k)}\left(\Xi_s-\Xi_{s-}\right)\mathrm{d}\tilde{J}_s^{(i)}
\end{aligned}
\tag{56}
$$

with the martingale part $\tilde{J}_s^{(i)}$ of $J_s^{(i)}$. Here, we use the notation $\Xi_s=\Xi\left(s,\left\{X_s^{(i)}\right\}_{i=1,2,3,\ldots,n}\right)$ and $\Xi_{s-}$ represents its left limit, and the summation $\sum_{s\in(t,t+k)}$ inside the other summation $\sum_{i=1}^{n}$ is with respect to all jump times of $J_s^{(i)}$. Then, we obtain

$$
\begin{aligned}
&\mathbb{E}_{\mathbb{Q}(\phi,\varphi)}\left[\Xi\left(t+k,\left\{X_{t+k}^{(i)}\right\}_{i=1,2,3,\ldots,n}\right)\right]-\Xi\left(t,\left\{X_{t}^{(i)}\right\}_{i=1,2,3,\ldots,n}\right) \\
&=\mathbb{E}_{\mathbb{Q}(\phi,\varphi)}\left[\int_t^{t+k}\left(\begin{array}{l}\dfrac{\partial \Xi_s}{\partial s}+\displaystyle\sum_{i=1}^{n}\left(\left(ac_i-r_iX_s^{(i)}+\sigma\sqrt{r_iX_s^{(i)}}\phi_s^{(i)}\right)\dfrac{\partial \Xi_s}{\partial X_s^{(i)}}+\dfrac{\sigma^2}{2}r_iX_s^{(i)}\dfrac{\partial^2 \Xi_s}{\partial X_s^{(i)}\partial X_s^{(i)}}\right) \\ +\displaystyle\sum_{i=1}^{n}c_i\int_0^{+\infty}\varphi_s^{(i)}\left(z_i\right)\left\{\Xi\left(s,\left\{X_s^{(j)}+\Delta_{i,j}z_j\right\}_{j=1,2,3,\ldots,n}\right)-\Xi_s\right\}v\left(\mathrm{d}z_i\right)\end{array}\right)\mathrm{d}s\right]
\end{aligned}
\tag{57}
$$

and hence

$$
\begin{aligned}
&\lim_{k\to+0}\frac{1}{k}\left\{\mathbb{E}_{\mathbb{Q}(\phi,\varphi)}\left[\Phi\left(\frac{\bar{\Psi}_{t+k}}{\bar{\Psi}_t}\right)\middle|\mathcal{F}_t\right]-\Phi(1)\right\} \\
&=\lim_{k\to+0}\frac{1}{k}\left\{\mathbb{E}_{\mathbb{Q}(\phi,\varphi)}\left[\Xi\left(t+k,\left\{X_{t+k}^{(i)}\right\}_{i=1,2,3,\ldots,n}\right)\right]-\Xi\left(t,\left\{X_{t}^{(i)}\right\}_{i=1,2,3,\ldots,n}\right)\right\} \\
&=\lim_{\left(s,\{y_i\}_{i=1,2,3,\ldots,n}\right)\to\left(t,\{x_i\}_{i=1,2,3,\ldots,n}\right)}\left(\begin{array}{l}\dfrac{\partial \hat{\Xi}}{\partial t}+\displaystyle\sum_{i=1}^{n}\left(\left(ac_i-r_ix_i+\sigma\sqrt{r_ix_i}\phi^{(i)}\right)\dfrac{\partial \hat{\Xi}}{\partial x_i}+\dfrac{\sigma^2}{2}r_ix_i\dfrac{\partial^2 \hat{\Xi}}{\partial x_i^2}\right) \\ +\displaystyle\sum_{i=1}^{n}c_i\int_0^{+\infty}\varphi^{(i)}\left(z_i\right)\left\{\hat{\Xi}\left(t,\left\{x_j+\Delta_{i,j}z_j\right\}_{j=1,2,3,\ldots,n}\right)-\hat{\Xi}\right\}v\left(\mathrm{d}z_i\right)\end{array}\right)
\end{aligned}
\tag{58}
$$

with $\hat{\Xi}=\Phi\left(\dfrac{\bar{\Psi}\left(t,\{x_i\}_{i=1,2,3,\ldots,n}\right)}{\bar{\Psi}\left(s,\{y_i\}_{i=1,2,3,\ldots,n}\right)}\right)$ and $\hat{\Xi}\left(t,\left\{x_j+\Delta_{i,j}z_j\right\}_{j=1,2,3,\ldots,n}\right)=\Phi\left(\dfrac{\bar{\Psi}\left(t,\left\{x_j+\Delta_{i,j}z_j\right\}_{j=1,2,3,\ldots,n}\right)}{\bar{\Psi}\left(s,\left\{y_j\right\}_{j=1,2,3,\ldots,n}\right)}\right)$.

By (54) and (58), we obtain

$$\sup_{\mathbb{A}(\phi,\varphi)}\left\{\begin{aligned}&\lim_{\left(s,\{y_i\}_{i=1,2,3,...,n}\right)\to\left(t,\{x_i\}_{i=1,2,3,...,n}\right)}\left(\begin{aligned}&\frac{\partial\hat{\Xi}}{\partial t}+\sum_{i=1}^{n}\left(\left(ac_i-r_ix_i+\sigma\sqrt{r_ix_i}\phi^{(i)}\right)\frac{\partial\hat{\Xi}}{\partial x_i}+\frac{\sigma^2}{2}r_ix_i\frac{\partial^2\hat{\Xi}}{\partial x_i^2}\right)\\&+\sum_{i=1}^{n}c_i\int_0^{+\infty}\varphi^{(i)}(z_i)\left\{\hat{\Xi}\left(t,\left\{x_j+\Delta_{i,j}z_j\right\}_{j=1,2,3,...,n}\right)-\hat{\Xi}\right\}\nu(\mathrm{d}z_i)\end{aligned}\right)\\&-\frac{1}{\lambda_\phi}\frac{1}{2}\sum_{i=1}^{n}\left(\phi^{(i)}\right)^2\\&-\frac{1}{\lambda_\varphi}\sum_{i=1}^{n}\int_0^{+\infty}\frac{1}{1-q}\left(1-q-\left(\varphi^{(i)}(z)\right)^q+q\varphi^{(i)}(z)\right)\nu(\mathrm{d}z)\end{aligned}\right\}=0 \quad (59)$$

with the set $\mathbb{A}(\phi,\varphi)$ given by

$$\mathbb{A}(\phi,\varphi)=\left\{(\phi,\varphi)\middle|\phi=\left\{\phi^{(i)}\right\}_{i=1,2,3,...,n}\in\mathbb{R}^n,\varphi=\left\{\varphi^{(i)}\right\}_{i=1,2,3,...,n}:(0,+\infty)^n\to(0,+\infty)^n\right\}. \quad (60)$$

We formally have

$$\begin{aligned}&\lim_{\left(s,\{y_i\}_{i=1,2,3,...,n}\right)\to\left(t,\{x_i\}_{i=1,2,3,...,n}\right)}\frac{\partial\hat{\Xi}}{\partial t}\\&=\lim_{\left(s,\{y_i\}_{i=1,2,3,...,n}\right)\to\left(t,\{x_i\}_{i=1,2,3,...,n}\right)}\frac{\partial}{\partial t}\Phi\left(\frac{\bar{\Psi}\left(t,\{x_i\}_{i=1,2,3,...,n}\right)}{\bar{\Psi}\left(s,\{y_i\}_{i=1,2,3,...,n}\right)}\right)\\&=\lim_{\left(s,\{y_i\}_{i=1,2,3,...,n}\right)\to\left(t,\{x_i\}_{i=1,2,3,...,n}\right)}\Phi'\left(\frac{\bar{\Psi}\left(t,\{x_i\}_{i=1,2,3,...,n}\right)}{\bar{\Psi}\left(s,\{y_i\}_{i=1,2,3,...,n}\right)}\right)\frac{1}{\bar{\Psi}\left(s,\{y_i\}_{i=1,2,3,...,n}\right)}\frac{\partial\bar{\Psi}\left(t,\{x_i\}_{i=1,2,3,...,n}\right)}{\partial t},\\&=\Phi'(1)\frac{1}{\bar{\Psi}\left(t,\{x_i\}_{i=1,2,3,...,n}\right)}\frac{\partial\bar{\Psi}\left(t,\{x_i\}_{i=1,2,3,...,n}\right)}{\partial t}\end{aligned} \quad (61)$$

$$\begin{aligned}&\lim_{\left(s,\{y_i\}_{i=1,2,3,...,n}\right)\to\left(t,\{x_i\}_{i=1,2,3,...,n}\right)}\frac{\partial\hat{\Xi}}{\partial x_j}\\&=\lim_{\left(s,\{y_i\}_{i=1,2,3,...,n}\right)\to\left(t,\{x_i\}_{i=1,2,3,...,n}\right)}\frac{\partial}{\partial x_j}\Phi\left(\frac{\bar{\Psi}\left(t,\{x_i\}_{i=1,2,3,...,n}\right)}{\bar{\Psi}\left(s,\{y_i\}_{i=1,2,3,...,n}\right)}\right)\\&=\lim_{\left(s,\{y_i\}_{i=1,2,3,...,n}\right)\to\left(t,\{x_i\}_{i=1,2,3,...,n}\right)}\Phi'\left(\frac{\bar{\Psi}\left(t,\{x_i\}_{i=1,2,3,...,n}\right)}{\bar{\Psi}\left(s,\{y_i\}_{i=1,2,3,...,n}\right)}\right)\frac{1}{\bar{\Psi}\left(s,\{y_i\}_{i=1,2,3,...,n}\right)}\frac{\partial\bar{\Psi}\left(t,\{x_i\}_{i=1,2,3,...,n}\right)}{\partial x_j},\\&=\Phi'(1)\frac{1}{\bar{\Psi}\left(t,\{x_i\}_{i=1,2,3,...,n}\right)}\frac{\partial\bar{\Psi}\left(t,\{x_i\}_{i=1,2,3,...,n}\right)}{\partial x_j}\end{aligned} \quad (62)$$

$$\lim_{\left(s,\{y_i\}_{i=1,2,3,\ldots,n}\right)\to\left(t,\{x_i\}_{i=1,2,3,\ldots,n}\right)} \frac{\partial^2 \hat{\Xi}}{\partial x_j^2}$$

$$= \lim_{\left(s,\{y_i\}_{i=1,2,3,\ldots,n}\right)\to\left(t,\{x_i\}_{i=1,2,3,\ldots,n}\right)} \frac{\partial^2}{\partial x_j^2}\Phi\left(\frac{\bar{\Psi}\left(t,\{x_i\}_{i=1,2,3,\ldots,n}\right)}{\bar{\Psi}\left(s,\{y_i\}_{i=1,2,3,\ldots,n}\right)}\right)$$

$$= \lim_{\left(s,\{y_i\}_{i=1,2,3,\ldots,n}\right)\to\left(t,\{x_i\}_{i=1,2,3,\ldots,n}\right)} \frac{\partial}{\partial x_j}\left\{\Phi'\left(\frac{\bar{\Psi}\left(t,\{x_i\}_{i=1,2,3,\ldots,n}\right)}{\bar{\Psi}\left(s,\{y_i\}_{i=1,2,3,\ldots,n}\right)}\right)\frac{1}{\bar{\Psi}\left(s,\{y_i\}_{i=1,2,3,\ldots,n}\right)}\frac{\partial\bar{\Psi}\left(t,\{x_i\}_{i=1,2,3,\ldots,n}\right)}{\partial x_j}\right\}$$

$$= \lim_{\left(s,\{y_i\}_{i=1,2,3,\ldots,n}\right)\to\left(t,\{x_i\}_{i=1,2,3,\ldots,n}\right)} \left\{\begin{array}{l}\Phi''\left(\frac{\bar{\Psi}\left(t,\{x_i\}_{i=1,2,3,\ldots,n}\right)}{\bar{\Psi}\left(s,\{y_i\}_{i=1,2,3,\ldots,n}\right)}\right)\left(\frac{1}{\bar{\Psi}\left(s,\{y_i\}_{i=1,2,3,\ldots,n}\right)}\right)^2\left(\frac{\partial\bar{\Psi}\left(t,\{x_i\}_{i=1,2,3,\ldots,n}\right)}{\partial x_j}\right)^2 \\ +\Phi'\left(\frac{\bar{\Psi}\left(t,\{x_i\}_{i=1,2,3,\ldots,n}\right)}{\bar{\Psi}\left(s,\{y_i\}_{i=1,2,3,\ldots,n}\right)}\right)\frac{1}{\bar{\Psi}\left(s,\{y_i\}_{i=1,2,3,\ldots,n}\right)}\frac{\partial^2\bar{\Psi}\left(t,\{x_i\}_{i=1,2,3,\ldots,n}\right)}{\partial x_j^2}\end{array}\right\}$$

$$= \Phi''(1)\left(\frac{1}{\bar{\Psi}\left(t,\{x_i\}_{i=1,2,3,\ldots,n}\right)}\right)^2\left(\frac{\partial\bar{\Psi}\left(t,\{x_i\}_{i=1,2,3,\ldots,n}\right)}{\partial x_j}\right)^2 + \Phi'(1)\frac{1}{\bar{\Psi}\left(t,\{x_i\}_{i=1,2,3,\ldots,n}\right)}\frac{\partial^2\bar{\Psi}\left(t,\{x_i\}_{i=1,2,3,\ldots,n}\right)}{\partial x_j^2}, \quad (63)$$

and

$$\lim_{\left(s,\{y_i\}_{i=1,2,3,\ldots,n}\right)\to\left(t,\{x_i\}_{i=1,2,3,\ldots,n}\right)} \left\{\hat{\Xi}\left(t,\{x_j+\Delta_{i,j}z_j\}_{j=1,2,3,\ldots,n}\right)-\hat{\Xi}\right\}$$

$$= \lim_{\left(s,\{y_i\}_{i=1,2,3,\ldots,n}\right)\to\left(t,\{x_i\}_{i=1,2,3,\ldots,n}\right)} \left\{\Phi\left(\frac{\bar{\Psi}\left(t,\{x_j+\Delta_{i,j}z_j\}_{j=1,2,3,\ldots,n}\right)}{\bar{\Psi}\left(s,\{y_i\}_{i=1,2,3,\ldots,n}\right)}\right)-\Phi\left(\frac{\bar{\Psi}\left(t,\{x_i\}_{i=1,2,3,\ldots,n}\right)}{\bar{\Psi}\left(s,\{y_i\}_{i=1,2,3,\ldots,n}\right)}\right)\right\}$$

$$= \Phi\left(\frac{\bar{\Psi}\left(t,\{x_j+\Delta_{i,j}z_j\}_{j=1,2,3,\ldots,n}\right)}{\bar{\Psi}\left(t,\{x_i\}_{i=1,2,3,\ldots,n}\right)}\right)-\Phi(1) \quad . \quad (64)$$

$$= \Phi\left(\frac{\bar{\Psi}\left(t,\{x_j+\Delta_{i,j}z_j\}_{j=1,2,3,\ldots,n}\right)}{\bar{\Psi}\left(t,\{x_i\}_{i=1,2,3,\ldots,n}\right)}\right)-1$$

By (61)–(64), (59) is rewritten as follows:

$$\sup_{\mathbb{A}(\phi,\varphi)}\left\{\begin{aligned}&\frac{\Phi'(1)}{\bar{\Psi}}\left(\frac{\partial\bar{\Psi}}{\partial t}+\sum_{i=1}^{n}\left(\left(ac_i-r_ix_i+\sigma\sqrt{r_ix_i}\phi^{(i)}\right)\frac{\partial\bar{\Psi}}{\partial x_i}+\frac{\sigma^2}{2}r_ix_i\frac{\partial^2\bar{\Psi}}{\partial x_i^2}\right)\right)\\&+\frac{\Phi''(1)}{\bar{\Psi}^2}\sum_{i=1}^{n}\frac{\sigma^2}{2}r_ix_i\left(\frac{\partial\bar{\Psi}}{\partial x_i}\right)^2\\&+\sum_{i=1}^{n}c_i\int_0^{+\infty}\varphi^{(i)}(z_i)\left\{\Phi\left(\frac{\bar{\Psi}\left(t,\left\{x_j+\Delta_{i,j}z_j\right\}_{j=1,2,3,\ldots,n}\right)}{\bar{\Psi}\left(t,\left\{x_j\right\}_{j=1,2,3,\ldots,n}\right)}\right)-1\right\}v(\mathrm{d}z_i)\\&-\frac{1}{\lambda_\phi}\frac{1}{2}\sum_{i=1}^{n}\left(\phi^{(i)}\right)^2\\&-\frac{1}{\lambda_\varphi}\sum_{i=1}^{n}\int_0^{+\infty}\frac{1}{1-q}\left(1-q-\left(\varphi^{(i)}(z)\right)^q+q\varphi^{(i)}(z)\right)v(\mathrm{d}z)\end{aligned}\right\}=0 \tag{65}$$

with the notation $\bar{\Psi}=\bar{\Psi}\left(t,\left\{x_i\right\}_{i=1,2,3,\ldots,n}\right)$. This equation becomes

$$\begin{aligned}&\frac{\Phi'(1)}{\bar{\Psi}}\left(\frac{\partial\bar{\Psi}}{\partial t}+\sum_{i=1}^{n}\left(\left(ac_i-r_ix_i\right)\frac{\partial\bar{\Psi}}{\partial x_i}+\frac{\sigma^2}{2}r_ix_i\frac{\partial^2\bar{\Psi}}{\partial x_i^2}\right)\right)+\frac{\Phi''(1)}{\bar{\Psi}^2}\sum_{i=1}^{n}\frac{\sigma^2}{2}r_ix_i\left(\frac{\partial\bar{\Psi}}{\partial x_i}\right)^2\\&+\sup_{\mathbb{A}(\phi,\varphi)}\left\{\frac{\Phi'(1)}{\bar{\Psi}}\sum_{i=1}^{n}\sigma\sqrt{r_ix_i}\phi^{(i)}\frac{\partial\bar{\Psi}}{\partial x_i}-\frac{1}{\lambda_\phi}\frac{1}{2}\sum_{i=1}^{n}\left(\phi^{(i)}\right)^2\right\}\\&+\sup_{\mathbb{A}(\phi,\varphi)}\left\{\begin{aligned}&\sum_{i=1}^{n}c_i\int_0^{+\infty}\varphi^{(i)}(z_i)\left\{\Phi\left(\frac{\bar{\Psi}\left(t,\left\{x_j+\Delta_{i,j}z_j\right\}_{j=1,2,3,\ldots,n}\right)}{\bar{\Psi}\left(t,\left\{x_j\right\}_{j=1,2,3,\ldots,n}\right)}\right)-1\right\}v(\mathrm{d}z_i)\\&-\frac{1}{\lambda_\varphi}\sum_{i=1}^{n}\int_0^{+\infty}\frac{1}{1-q}\left(1-q-\left(\varphi^{(i)}(z)\right)^q+q\varphi^{(i)}(z)\right)v(\mathrm{d}z)\end{aligned}\right\}\\&=0\end{aligned}\quad. \tag{66}$$

The second line of (66) is rewritten as follows:

$$\begin{aligned}&\sup_{\mathbb{A}(\phi,\varphi)}\left\{\frac{\Phi'(1)}{\bar{\Psi}}\sum_{i=1}^{n}\sigma\sqrt{r_ix_i}\phi^{(i)}\frac{\partial\bar{\Psi}}{\partial x_i}-\frac{1}{\lambda_\phi}\frac{1}{2}\sum_{i=1}^{n}\left(\phi^{(i)}\right)^2\right\}\\&=\frac{\Phi'(1)}{\bar{\Psi}}\sup_{\mathbb{A}(\phi,\varphi)}\left\{\sum_{i=1}^{n}\sigma\sqrt{r_ix_i}\phi^{(i)}\frac{\partial\bar{\Psi}}{\partial x_i}-\frac{\bar{\Psi}}{\lambda_\phi\Phi'(1)}\frac{1}{2}\sum_{i=1}^{n}\left(\phi^{(i)}\right)^2\right\}\\&=\frac{\Phi'(1)^2\lambda_\phi}{\bar{\Psi}^2}\frac{1}{2}\sum_{i=1}^{n}\sigma r_ix_i\left(\frac{\partial\bar{\Psi}}{\partial x_i}\right)^2\end{aligned} \tag{67}$$

with the maximizer

$$\left\{\phi^{(i),*}\right\}_{i=1,2,3,\ldots,n}=\left\{\frac{\lambda_\phi\Phi'(1)}{\bar{\Psi}}\sigma\sqrt{r_ix_i}\frac{\partial\bar{\Psi}}{\partial x_i}\right\}_{i=1,2,3,\ldots,n}. \tag{68}$$

The third line of (66) is rewritten as follows:

$$\sup_{\mathbb{A}(\phi,\varphi)}\left\{\begin{aligned}&\sum_{i=1}^{n}c_i\int_0^{+\infty}\varphi^{(i)}(z_i)\left\{\Phi\left(\frac{\bar{\Psi}\left(t,\left\{x_j+\Delta_{i,j}z_j\right\}_{j=1,2,3,\ldots,n}\right)}{\bar{\Psi}\left(t,\left\{x_j\right\}_{j=1,2,3,\ldots,n}\right)}\right)-1\right\}v(\mathrm{d}z_i)\\&-\frac{1}{\lambda_\varphi}\sum_{i=1}^{n}\int_0^{+\infty}\frac{1}{1-q}\left(1-q-\left(\varphi^{(i)}(z)\right)^q+q\varphi^{(i)}(z)\right)v(\mathrm{d}z)\end{aligned}\right\}$$
$$=\frac{1}{\lambda_\varphi}\sum_{i=1}^{n}c_i\int_0^{+\infty}\left(\exp_q\left\{\lambda_\varphi\Phi\left(\frac{\bar{\Psi}\left(t,\left\{x_j+\Delta_{i,j}z_j\right\}_{j=1,2,3,\ldots,n}\right)}{\bar{\Psi}\left(t,\left\{x_j\right\}_{j=1,2,3,\ldots,n}\right)}\right)-\lambda_\varphi\right\}-1\right)v(\mathrm{d}z_i) \tag{69}$$

with the maximizer

$$\left\{\varphi^{(i),*}(z_i)\right\}_{i=1,2,3,\ldots,n}=\left\{\exp_q\left\{\lambda_\varphi\Phi\left(\frac{\bar{\Psi}\left(t,\left\{x_j+\Delta_{i,j}z_j\right\}_{j=1,2,3,\ldots,n}\right)}{\bar{\Psi}\left(t,\left\{x_j\right\}_{j=1,2,3,\ldots,n}\right)}\right)-\lambda_\varphi\right\}\right\}_{i=1,2,3,\ldots,n}. \tag{70}$$

Now, the following HJB equation is obtained:

$$\frac{\partial\bar{\Psi}}{\partial t}+\sum_{i=1}^{n}\left(\left(ac_i-r_ix_i\right)\frac{\partial\bar{\Psi}}{\partial x_i}+\frac{\sigma^2}{2}r_ix_i\frac{\partial^2\bar{\Psi}}{\partial x_i^2}\right)+\frac{\Phi'(1)^2\lambda_\phi+\Phi''(1)}{\Phi'(1)}\frac{1}{\bar{\Psi}}\sum_{i=1}^{n}\frac{\sigma^2}{2}r_ix_i\left(\frac{\partial\bar{\Psi}}{\partial x_i}\right)^2$$
$$+\frac{1}{\lambda_\varphi\Phi'(1)}\bar{\Psi}\sum_{i=1}^{n}c_i\int_0^{+\infty}\left(\exp_q\left\{\lambda_\varphi\Phi\left(\frac{\bar{\Psi}\left(t,\left\{x_j+\Delta_{i,j}z_j\right\}_{j=1,2,3,\ldots,n}\right)}{\bar{\Psi}\left(t,\left\{x_j\right\}_{j=1,2,3,\ldots,n}\right)}\right)-\lambda_\varphi\right\}-1\right)v(\mathrm{d}z_i)=0. \tag{71}$$

This is the desired equation (29). Thus, the first part of the proof is complete.

The second part of the proof involves a substitution of (32). Here, we check the sufficiency of the conditions $0<q<1$ and $0<p<\min\left\{\frac{2}{\sigma^2},\beta\right\}$, as well as that of (28) with $\frac{mp}{1-q}<\beta$. Firstly, it follows from (89) that $0<\bar{\rho}_s^{(i)}\le p$ for all $s\le T$ by the assumption. Secondly, the condition $0<q<1$ combined with $0<\bar{\rho}_s^{(i)}\le p$ for all $s\le T$ ensures that the argument inside each "$\exp_q$" in (71) is nonnegative, and the equation is well-defined. Finally, for some sufficiently large constant $C>0$, we obtain:

$$
\begin{aligned}
\exp_q\left\{\lambda_\varphi \Phi\left(\frac{\bar{\Psi}\left(t,\left\{x_j+\Delta_{i,j}z_j\right\}_{j=1,2,3,\ldots,n}\right)}{\bar{\Psi}\left(t,\left\{x_j\right\}_{j=1,2,3,\ldots,n}\right)}\right)-\lambda_\varphi\right\} &= \exp_q\left\{\lambda_\varphi \Phi\left(e^{\bar{\rho}_t^{(i)}z}\right)-\lambda_\varphi\right\} \\
&= \left(1+(1-a)\lambda_\varphi\left(\Phi\left(e^{\bar{\rho}_t^{(i)}z}\right)-1\right)\right)^{\frac{1}{1-q}} \\
&\leq \left(1+(1-a)\lambda_\varphi\left(\Phi\left(e^{pz}\right)-1\right)\right)^{\frac{1}{1-q}} \\
&\leq C\Phi\left(e^{pz}\right)^{\frac{1}{1-q}} \\
&\leq C\left(C_0 e^{mpz}+C_1\right)^{\frac{1}{1-q}} \\
&\leq C_0 e^{\frac{mp}{1-q}z}+C_1
\end{aligned}
\tag{72}
$$

with the abuse of the notation for sufficiently large generic constants $C_0, C_1 > 0$. Subsequently, owing to the regularity of $\Phi$, (2) implies that the last integral is bounded under $m \in \left(1, \frac{\beta}{p}(1-q)\right)$. This proves the second step in the proof of **Proposition 1**:

***Remark A1*** For $q=1$, we formally obtain the HJB equation having the integral term

$$
\begin{aligned}
&\int_0^{+\infty}\left(\exp\left\{\lambda_\varphi \Phi\left(\frac{\bar{\Psi}\left(t,\left\{x_j+\Delta_{i,j}z_j\right\}_{j=1,2,3,\ldots,n}\right)}{\bar{\Psi}\left(t,\left\{x_j\right\}_{j=1,2,3,\ldots,n}\right)}\right)-\lambda_\varphi\right\}-1\right)v\left(\mathrm{d}z_i\right) \\
&= \int_0^{+\infty}\left(\exp\left\{\lambda_\varphi \Phi\left(e^{\bar{\rho}_t^{(i)}z_i}\right)-\lambda_\varphi\right\}-1\right)v\left(\mathrm{d}z_i\right) \\
&= \gamma\int_0^{+\infty} z^{-(1+\alpha)}\left(\exp\left\{\lambda_\varphi \Phi\left(e^{\bar{\rho}_t^{(i)}z_i}\right)-\lambda_\varphi-\beta z\right\}-1\right)\mathrm{d}z_i
\end{aligned}
, \tag{73}
$$

which is not integrable at any $t<T$ because $\exp\left\{\lambda_\varphi \Phi\left(e^{\bar{\rho}_t^{(i)}z_i}\right)-\lambda_\varphi-\beta z\right\}$ grows doubly exponentially for $z>0$. As a result, the case $q=1$ is not included in **Proposition 1**.

## Appendix B. Proof of Proposition 2

We assume $q>1$ to prove that $q=1$ remains essentially constant. The proof is essentially identical to that of **Proposition 1**—hence, only the key differences are explained here. Firstly, Orlicz risk deals with the infimum, not the supremum. The HJB reads

$$\begin{aligned}
&\frac{\Phi'(1)}{\underline{\Psi}}\left(\frac{\partial\underline{\Psi}}{\partial t}+\sum_{i=1}^{n}\left(\left(ac_i-r_ix_i\right)\frac{\partial\underline{\Psi}}{\partial x_i}+\frac{\sigma^2}{2}r_ix_i\frac{\partial^2\underline{\Psi}}{\partial x_i^2}\right)\right)+\Phi''(1)\frac{1}{\underline{\Psi}^2}\sum_{i=1}^{n}\frac{\sigma^2}{2}r_ix_i\left(\frac{\partial\underline{\Psi}}{\partial x_i}\right)^2\\
&+\inf_{\mathbb{A}(\phi,\varphi)}\left\{\frac{\Phi'(1)}{\underline{\Psi}}\sum_{i=1}^{n}\sigma\sqrt{r_ix_i}\phi^{(i)}\frac{\partial\underline{\Psi}}{\partial x_i}+\frac{1}{\lambda_\phi}\frac{1}{2}\sum_{i=1}^{n}\left(\phi^{(i)}\right)^2\right\}\\
&+\inf_{\mathbb{A}(\phi,\varphi)}\left\{\begin{array}{l}\displaystyle\sum_{i=1}^{n}c_i\int_0^{+\infty}\varphi^{(i)}(z_i)\left\{\Phi\left(\frac{\underline{\Psi}\left(t,\left\{x_j+\Delta_{i,j}z_j\right\}_{j=1,2,3,\ldots,n}\right)}{\underline{\Psi}\left(t,\left\{x_j\right\}_{j=1,2,3,\ldots,n}\right)}\right)-1\right\}v(\mathrm{d}z_i)\\ \displaystyle+\frac{1}{\lambda_\varphi}\sum_{i=1}^{n}\int_0^{+\infty}\frac{1}{1-q}\left(1-q-\left(\varphi^{(i)}(z)\right)^q+q\varphi^{(i)}(z)\right)v(\mathrm{d}z)\end{array}\right\}\\
&=0
\end{aligned} \tag{74}$$

We have

$$\inf_{\mathbb{A}(\phi,\varphi)}\left\{\frac{\Phi'(1)}{\underline{\Psi}}\sum_{i=1}^{n}\sigma\sqrt{r_ix_i}\phi^{(i)}\frac{\partial\underline{\Psi}}{\partial x_i}+\frac{1}{\lambda_\phi}\frac{1}{2}\sum_{i=1}^{n}\left(\phi^{(i)}\right)^2\right\}=-\frac{\Phi'(1)^2\lambda_\phi}{\underline{\Psi}^2}\frac{1}{2}\sum_{i=1}^{n}\sigma r_ix_i\left(\frac{\partial\underline{\Psi}}{\partial x_i}\right)^2 \tag{75}$$

with the minimizer

$$\left\{\phi^{(i),*}\right\}_{i=1,2,3,\ldots,n}=\left\{-\frac{\lambda_\phi\Phi'(1)}{\underline{\Psi}}\sigma\sqrt{r_ix_i}\frac{\partial\underline{\Psi}}{\partial x_i}\right\}_{i=1,2,3,\ldots,n} \tag{76}$$

as well as

$$\begin{aligned}
&\inf_{\mathbb{A}(\phi,\varphi)}\left\{\begin{array}{l}\displaystyle\sum_{i=1}^{n}c_i\int_0^{+\infty}\varphi^{(i)}(z_i)\left\{\Phi\left(\frac{\underline{\Psi}\left(t,\left\{x_j+\Delta_{i,j}z_j\right\}_{j=1,2,3,\ldots,n}\right)}{\underline{\Psi}\left(t,\left\{x_j\right\}_{j=1,2,3,\ldots,n}\right)}\right)-1\right\}v(\mathrm{d}z_i)\\ \displaystyle+\frac{1}{\lambda_\varphi}\sum_{i=1}^{n}\int_0^{+\infty}\frac{1}{1-q}\left(1-q-\left(\varphi^{(i)}(z)\right)^q+q\varphi^{(i)}(z)\right)v(\mathrm{d}z)\end{array}\right\}\\
&=-\frac{1}{\lambda_\varphi}\sum_{i=1}^{n}c_i\int_0^{+\infty}\left(\exp_q\left\{-\lambda_\varphi\Phi\left(\frac{\underline{\Psi}\left(t,\left\{x_j+\Delta_{i,j}z_j\right\}_{j=1,2,3,\ldots,n}\right)}{\underline{\Psi}\left(t,\left\{x_j\right\}_{j=1,2,3,\ldots,n}\right)}\right)+\lambda_\varphi\right\}-1\right)v(\mathrm{d}z_i)
\end{aligned} \tag{77}$$

with the minimizer

$$\left\{\varphi^{(i),*}(z_i)\right\}_{i=1,2,3,\ldots,n}=\left\{\exp_q\left\{-\lambda_\varphi\Phi\left(\frac{\underline{\Psi}\left(t,\left\{x_j+\Delta_{i,j}z_j\right\}_{j=1,2,3,\ldots,n}\right)}{\underline{\Psi}\left(t,\left\{x_j\right\}_{j=1,2,3,\ldots,n}\right)}\right)+\lambda_\varphi\right\}\right\}_{i=1,2,3,\ldots,n}. \tag{78}$$

The HJB equation (74) yields

$$\frac{\partial \underline{\Psi}}{\partial t}+\sum_{i=1}^{n}\left(\left(ac_i-r_ix_i\right)\frac{\partial \underline{\Psi}}{\partial x_i}+\frac{\sigma^2}{2}r_ix_i\frac{\partial^2 \underline{\Psi}}{\partial x_i^2}\right)-\frac{\Phi'(1)^2\lambda_\phi-\Phi''(1)}{\Phi'(1)}\frac{1}{\underline{\Psi}}\sum_{i=1}^{n}\frac{\sigma^2}{2}r_ix_i\left(\frac{\partial \underline{\Psi}}{\partial x_i}\right)^2$$
$$-\frac{1}{\lambda_\varphi\Phi'(1)}\underline{\Psi}\sum_{i=1}^{n}c_i\int_0^{+\infty}\left(\exp_q\left\{-\lambda_\varphi\Phi\left(\frac{\underline{\Psi}\left(t,\left\{x_j+\Delta_{i,j}z_j\right\}_{j=1,2,3,\dots,n}\right)}{\underline{\Psi}\left(t,\left\{x_j\right\}_{j=1,2,3,\dots,n}\right)}\right)+\lambda_\varphi\right\}-1\right)v(\mathrm{d}z_i)=0 \quad , \tag{79}$$

which is identical to (38).

All that remains to be verified is the well-defined nature of the last integral in (79) by the closed-form solution in (42) and (43). Let us assume $q>1$. Then, we have

$$\begin{aligned}\exp_q\left\{-\lambda_\varphi\Phi\left(\frac{\underline{\Psi}\left(t,\left\{x_j+\Delta_{i,j}z_j\right\}_{j=1,2,3,\dots,n}\right)}{\underline{\Psi}\left(t,\left\{x_j\right\}_{j=1,2,3,\dots,n}\right)}\right)+\lambda_\varphi\right\}&=\exp_q\left\{-\lambda_\varphi\Phi\left(e^{\underline{\rho}_t^{(i)}z_i}\right)+\lambda_\varphi\right\}\\&=\frac{1}{\left(1+(q-1)\lambda_\varphi\left(\Phi\left(e^{\underline{\rho}_t^{(i)}z_i}\right)-1\right)\right)^{\frac{1}{q-1}}}\,,\\&\in(0,1)\end{aligned} \tag{80}$$

where the argument of " $\exp_q$ " is well-defined because of $q>1$. For a small $z_i>0$, we have

$$\exp_q\left\{-\lambda_\varphi\Phi\left(\frac{\underline{\Psi}\left(t,\left\{x_j+\Delta_{i,j}z_j\right\}_{j=1,2,3,\dots,n}\right)}{\underline{\Psi}\left(t,\left\{x_j\right\}_{j=1,2,3,\dots,n}\right)}\right)+\lambda_\varphi\right\}-1=O(z_i)\,, \tag{81}$$

which, combined with (80) implies that the integrand in the last integral of (79) is integrable. Hence, the integral is well defined.

***Remark B1*** We present the derivation process for SDEs (51) and (52). In the upper bound case, the worst-case uncertainty in a stationary state can be obtained from (68) as follows:

$$\phi_t^{(i),*}=\frac{\lambda_\phi\Phi'(1)}{\underline{\Psi}}\sigma\sqrt{r_iX_t^{(i)}}\frac{\partial\underline{\Psi}}{\partial X_t^{(i)}}=\lambda_\phi\Phi'(1)\sigma\underline{\rho}_t^{(i)}\sqrt{r_iX_t^{(i)}}\,. \tag{82}$$

Similarly, we have

$$\varphi_t^{(i),*}(z_i)=\exp_q\left\{-\lambda_\varphi\Phi\left(e^{\underline{\rho}_t^{(i)}z}\right)+\lambda_\varphi\right\}. \tag{83}$$

As we are interested in a stationary state, we assume $t=T$ and take $T$ to be sufficiently large, and hence arrive at $\underline{\rho}_t^{(i)}=p$ in (82) and (83). Then, we obtain the SDE in the stationary state as follows:

$$\begin{aligned}\mathrm{d}X_t^{(i)}&=\left(ac_i-r_iX_t^{(i)}+\sigma\sqrt{r_iX_t^{(i)}}\phi_t^{(i),*}\right)\mathrm{d}t+\sigma\sqrt{r_iX_t^{(i)}}\mathrm{d}W_t^{(i)}+\mathrm{d}\underline{J}_t^{(i)}\\&=\left(ac_i-r_i(1+\xi)X_t^{(i)}\right)\mathrm{d}t+\sigma\sqrt{r_iX_t^{(i)}}\mathrm{d}W_t^{(i)}+\mathrm{d}\underline{J}_t^{(i)}\end{aligned}\quad , \tag{84}$$

which is identical to (51). The SDE in the lower bound case can be derived similarly. The moments and autocorrelation functions of the supJCIR process based on distorted SDEs can be obtained analytically, as in the case involving no uncertainty. For example, given $\xi \in (0,1)$, the autocorrelation function for the upper bound case (under the infinite-dimensional limit $n \to +\infty$) is

$$\mathrm{ACF}(h) = \frac{1}{\left(1+\theta(1-\xi)h\right)^{\omega-1}}, \quad h \geq 0. \tag{85}$$

Similarly, that for the lower bound case is

$$\mathrm{ACF}(h) = \frac{1}{\left(1+\theta(1+\xi)h\right)^{\omega-1}}, \quad h \geq 0. \tag{86}$$

**Appendix C. The infinite-dimensional limit**

The infinite-dimensional model consistent with the finite-dimensional one analyzed in the main text is explained here. Under the infinite-dimensional limit ( $n \to +\infty$ ), the SDE (8) is interpreted as an integration of a continuum of JCIR processes, and the corresponding HJB equations in **Propositions 1–2** as partial differential equations in infinite dimensions. This formulation relies on the formulation based on measure-valued processes [19,27,74].

The SDE system (8)–(9) is formally interpreted as

$$Y_t = \int_0^{+\infty} X_t(\mathrm{d}r),\quad t > 0 \tag{87}$$

with the Itô's SDE

$$\mathrm{d}X_t(\mathrm{d}r) = \left(a\pi(\mathrm{d}r) - rX_t(\mathrm{d}r)\right)\mathrm{d}t + \sigma\sqrt{rX_t(\mathrm{d}r)}\mathrm{d}B_t(\mathrm{d}r) + \mathrm{d}L_t(\mathrm{d}r),\quad r > 0, \tag{88}$$

where $\left(B_t(\mathrm{d}r)\right)_{t\geq 0}$ for $r > 0$ is a Brownian motion that is mutually independent for different $r = r_1, r_2 > 0$ and satisfies the product rule $\mathrm{d}B_t(\mathrm{d}r)\mathrm{d}B_t(\mathrm{d}r) = \pi(\mathrm{d}r)\mathrm{d}t$ . Moreover, $\left(L_t(\mathrm{d}r)\right)_{t\geq 0}$ is the space-time having Lévy process with the compensator $\pi(\mathrm{d}r)v(\mathrm{d}z)\mathrm{d}t$ ; i.e., Lévy bases [27]. The formulation (88) implies that each $X_t(\mathrm{d}r)$ has a magnitude of $\pi(\mathrm{d}r)$ , and their integration (87) is the order of 1. This kind of infinite-dimensional formulation is still under study [30].

The HJB equations corresponding to those of **Propositions 1–2** in the infinite-dimensional framework are of forms like those in the literature [30] and are not presented here because of their length. The important conclusion in the context of this study is that the infinite-dimensional counterparts of the closed-form solutions are integro-differential equations, which are easier to understand and handle. In this context, for mappings $\left(\bar{\tau}_t, \bar{\rho}_t(\cdot)\right)_{t\leq T}$ , with the abuse of notations we arrive at the representations

$$\bar{\rho}_t(r) = \frac{1}{\frac{\sigma^2\left(1+\bar{\lambda}\right)}{2} + \left(\frac{1}{p} - \frac{\sigma^2\left(1+\bar{\lambda}\right)}{2}\right)e^{r(T-t)}},\quad r > 0,\quad t \leq T \tag{89}$$

and

$$\frac{\mathrm{d}\bar{\tau}_t}{\mathrm{d}t} + a\int_0^{+\infty}\bar{\rho}_t(r)\pi(\mathrm{d}r) + \frac{1}{\lambda_\varphi\Phi'(1)}\int_0^{+\infty}\int_0^{+\infty}\left\{\exp_q\left(\lambda_\varphi\left(\Phi\left(e^{\bar{\rho}_t(r)z}\right)-1\right)\right)-1\right\}v(\mathrm{d}z)\pi(\mathrm{d}r) = 0,\quad t \leq T. \tag{90}$$

Then, for a generic sufficiently regular nonnegative mapping $y:(0,+\infty)\to(0,+\infty)$ , the exponential disutility becomes

$$\bar{F} = \exp\left(\int_0^{+\infty}\bar{\rho}_t(r)y(r)\mathrm{d}r + \bar{\tau}_t\right). \tag{91}$$

Similarly, for mapping $\left(\underline{\tau}_t, \underline{\rho}_t(\cdot)\right)_{t\leq T}$ , (42) and (43) become, with an abuse of notations:

$$\underline{\rho}_t(r) = \frac{1}{\frac{\sigma^2\left(1-\underline{\lambda}\right)}{2} + \left(\frac{1}{p} - \frac{\sigma^2\left(1-\underline{\lambda}\right)}{2}\right)e^{r(T-t)}},\quad r > 0,\quad t \leq T \tag{92}$$

and

$$\frac{\mathrm{d}\underline{\tau}_t}{\mathrm{d}t}+a\int_0^{+\infty}\underline{\rho}_t(r)\pi(\mathrm{d}r)-\frac{1}{\lambda_\varphi\Phi'(1)}\int_0^{+\infty}\int_0^{+\infty}\left\{\exp_q\left(-\lambda_\varphi\left(\Phi\left(e^{\underline{\rho}_t(r)z}\right)-1\right)\right)-1\right\}v(\mathrm{d}z)\pi(\mathrm{d}r)=0\,,\quad t\leq T\,. \tag{93}$$

For a generic sufficiently regular mapping $y:(0,+\infty)\to(0,+\infty)$, the exponential disutility becomes

$$\underline{F}=\exp\left(\int_0^{+\infty}\underline{\rho}_t(r)y(r)\mathrm{d}r+\underline{\tau}_t\right). \tag{94}$$

From this perspective, the closed-form solutions of **Propositions 1–2** are discretization of (91) and (94).

***Remark C1***: Formulation (90) can be simplified to a stationary state. We have:

$$\overline{\tau}_{-\infty}=\int_{-\infty}^{T}\left(a\int_0^{+\infty}\overline{\rho}_s(r)\pi(\mathrm{d}r)+\frac{1}{\lambda_\varphi\Phi'(1)}\int_0^{+\infty}\int_0^{+\infty}\left\{\exp_q\left(\lambda_\varphi\left(\Phi\left(e^{\overline{\rho}_s(r)z}\right)-1\right)\right)-1\right\}v(\mathrm{d}z)\pi(\mathrm{d}r)\right)\mathrm{d}s\,. \tag{95}$$

For the first term in (95), we have

$$\begin{aligned}\int_{-\infty}^{T}\int_0^{+\infty}\overline{\rho}_s(r)\pi(\mathrm{d}r)\mathrm{d}s&=\int_{-\infty}^{T}\int_0^{+\infty}\frac{1}{r}\overline{\rho}_{s/r}(r)\pi(\mathrm{d}r)\mathrm{d}s\\&=\int_{-\infty}^{T}\int_0^{+\infty}\frac{1}{r}\overline{\rho}_s(1)\pi(\mathrm{d}r)\mathrm{d}s\\&=\left(\int_0^{+\infty}\frac{1}{r}\pi(\mathrm{d}r)\right)\int_{-\infty}^{T}\overline{\rho}_s(1)\mathrm{d}s\\&=R\int_{-\infty}^{T}\overline{\rho}_s(1)\mathrm{d}s\end{aligned}\,. \tag{96}$$

Similarly, for the second term in (95), we have

$$\begin{aligned}&\int_{-\infty}^{T}\left(\frac{1}{\lambda_\varphi\Phi'(1)}\int_0^{+\infty}\int_0^{+\infty}\left\{\exp_q\left(\lambda_\varphi\left(\Phi\left(e^{\overline{\rho}_s(r)z}\right)-1\right)\right)-1\right\}v(\mathrm{d}z)\pi(\mathrm{d}r)\right)\mathrm{d}s\\&=\int_{-\infty}^{T}\left(\frac{1}{\lambda_\varphi\Phi'(1)}\int_0^{+\infty}\frac{1}{r}\int_0^{+\infty}\left\{\exp_q\left(\lambda_\varphi\left(\Phi\left(e^{\overline{\rho}_{s/r}(r)z}\right)-1\right)\right)-1\right\}v(\mathrm{d}z)\pi(\mathrm{d}r)\right)\mathrm{d}s\\&=\frac{R}{\lambda_\varphi\Phi'(1)}\int_{-\infty}^{T}\int_0^{+\infty}\left\{\exp_q\left(\lambda_\varphi\left(\Phi\left(e^{\overline{\rho}_s(1)z}\right)-1\right)\right)-1\right\}v(\mathrm{d}z)\mathrm{d}s\end{aligned}\,. \tag{97}$$

We therefore obtain:

$$\overline{\tau}_{-\infty}=R\left\{\int_{-\infty}^{T}\overline{\rho}_s(1)\mathrm{d}s+\frac{1}{\lambda_\varphi\Phi'(1)}\int_{-\infty}^{T}\int_0^{+\infty}\left\{\exp_q\left(\lambda_\varphi\left(\Phi\left(e^{\overline{\rho}_s(1)z}\right)-1\right)\right)-1\right\}v(\mathrm{d}z)\mathrm{d}s\right\}. \tag{98}$$

Similarly, we obtain:

$$\underline{\tau}_{-\infty}=R\left\{\int_{-\infty}^{T}\underline{\rho}_s(1)\mathrm{d}s-\frac{1}{\lambda_\varphi\Phi'(1)}\int_{-\infty}^{T}\int_0^{+\infty}\left\{\exp_q\left(-\lambda_\varphi\left(\Phi\left(e^{\underline{\rho}_s(1)z}\right)-1\right)\right)-1\right\}v(\mathrm{d}z)\mathrm{d}s\right\}. \tag{99}$$

Therefore, integration with respect to $r$ can be completed without any computational errors. These formulae are described in **Section 4**.

**Appendix D. Auxiliary computational results**

**D.1 On the model for $SO_4^{2-}$**

**Figure D1** depicts the upper (red, $\Phi = x^{3/2}$ ) and lower (blue, $\Phi = x^{1/2}$ ) bounds of the exponential disutility for $SO_4^{2-}$ with $p = 0.01$. Here, we use the identified model without considering the last term in equation (48). The weight $\lambda_\phi$ is taken to be $\lambda_\phi = 0.01 + 30k$ ( $k = 0,1,2,...,200$ ) for the upper bound and $\lambda_\phi = 0.01 + 100000k$ ( $k = 0,1,2,...,200$ ) for the lower bound. The upper and lower bounds behave similarly to those of the TN, as discussed in **Section 4** of the main text. Its dependence on the weight $\lambda_\phi$ is sharper for $SO_4^{2-}$ than that for TN because of its larger $\sigma$.

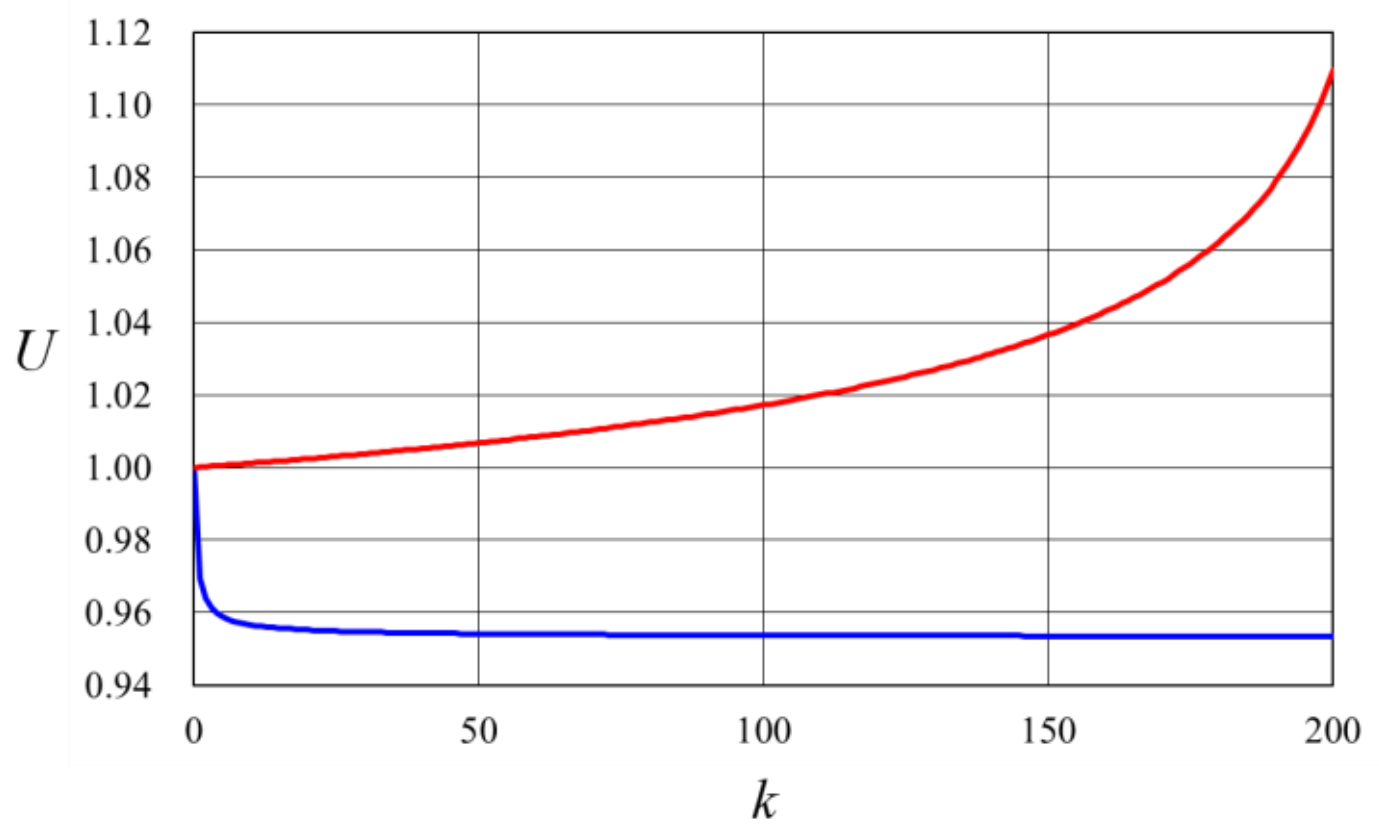

**Figure D1.** The upper bound (red) and lower bound (blue) of the exponential disutility for $SO_4^{2-}$ with $p = 0.01$. The curves are plotted as a function of $k$, where we have chosen $\lambda_\phi = 0.01 + 30k$ ( $k = 0,1,2,...,200$ ) for the upper bound and $\lambda_\phi = 0.01 + 100000k$ ( $k = 0,1,2,...,200$ ) for the lower bound. Here, the exponential utility is 1.04912 in the absence of model uncertainty.

**D.2 Moments under model uncertainties**

**Figures D2-D3** depict the computed averages and variances of TN in the upper and lower bound cases, respectively, where $p = 0.01$. The normalized average, $A$, and normalized variance, $V$, of TN in these figures are defined as follows:

$$A = \frac{\left.\text{Average}\right|_{\text{With uncertainty}}}{\left.\text{Average}\right|_{\text{No uncertainty}}} \text{ and } V = \frac{\left.\text{Variance}\right|_{\text{With uncertainty}}}{\left.\text{Variance}\right|_{\text{No uncertainty}}} . \tag{100}$$

These figures indicate that the normalized average and variance monotonically increase and decrease, respectively, for the ambiguity aversion parameters in the upper and lower bound cases.

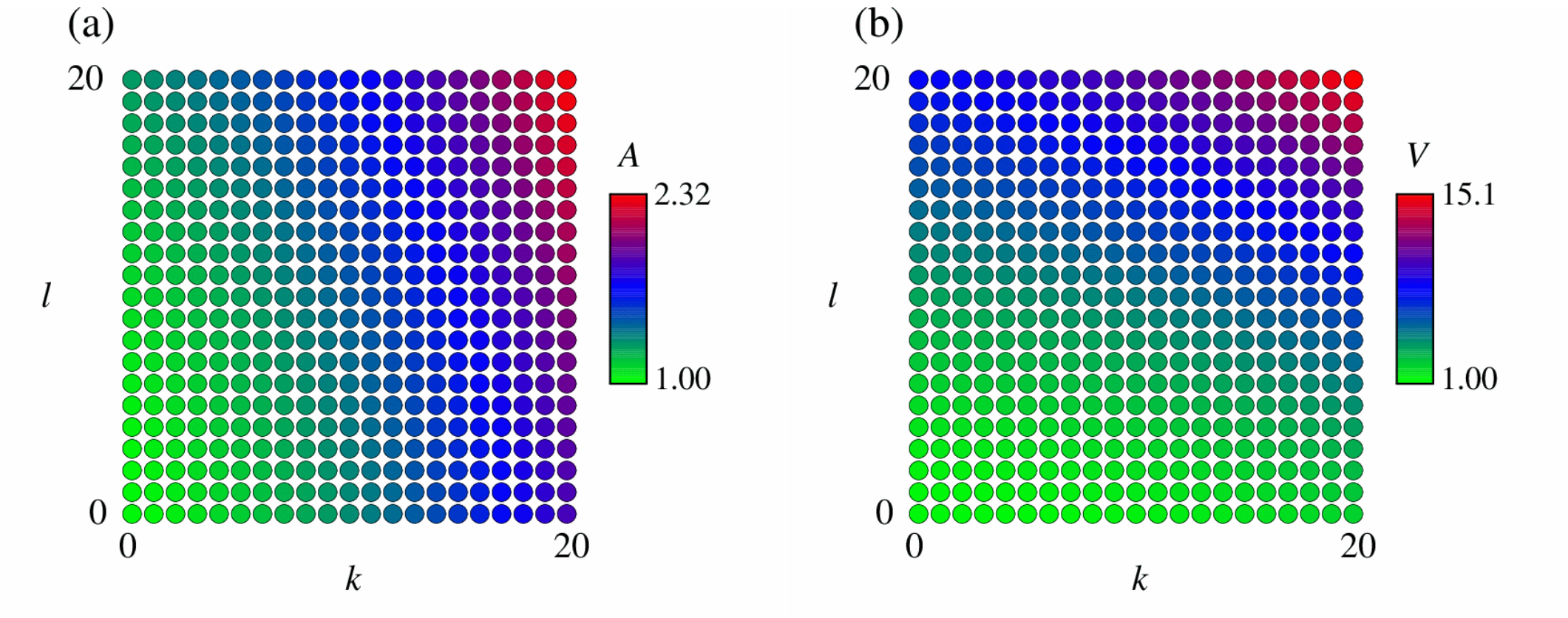

**Figure D2.** The computed (a) normalized average $A$, and (b) normalized variance $V$ in the upper bound, where $\Phi = x^{3/2}$, $q = 0.5$, and $(\lambda_\phi, \lambda_\varphi) = (0.01 + 5000k, 0.01 + 200l)$ ($k, l = 0,1,2,...,20$).

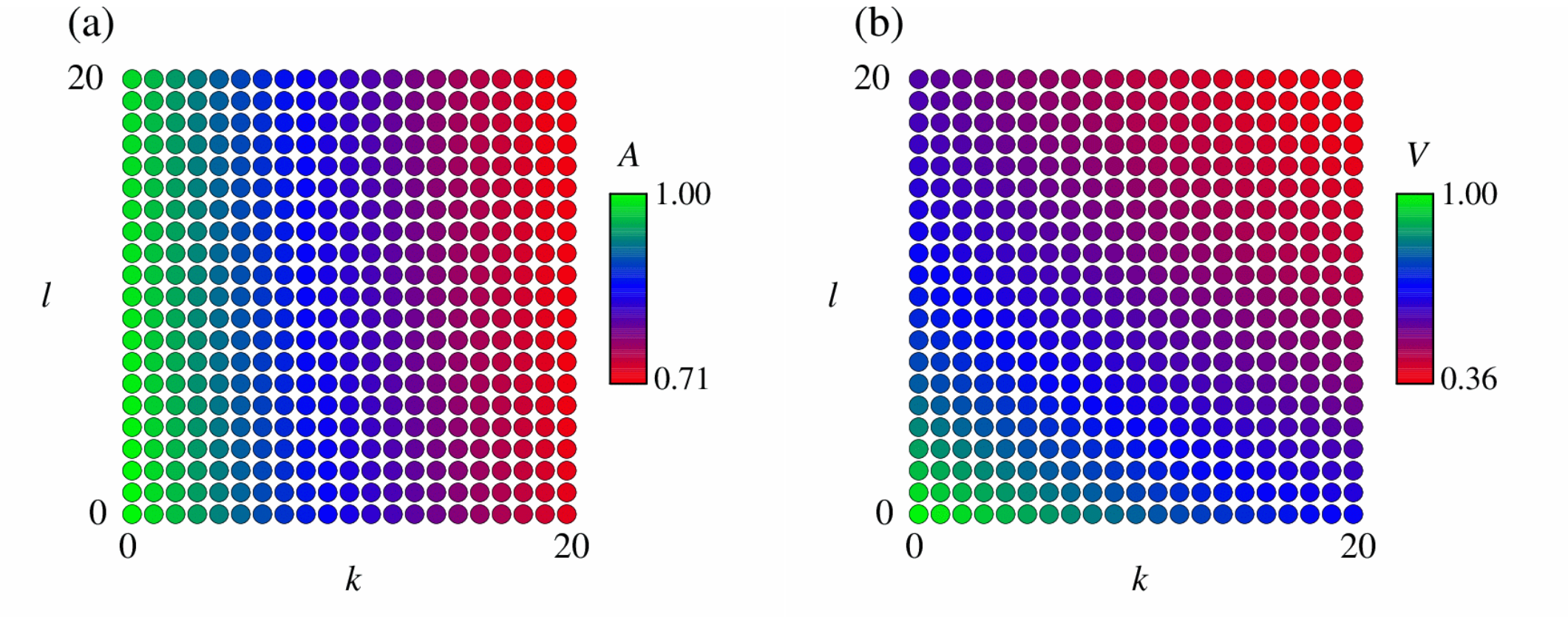

**Figure D3**. The computed (a) normalized average $A$ and (b) normalized variance $V$ in the lower bound case, where $\Phi = x^{3/2}$, $q = 1.5$, and $(\lambda_\phi, \lambda_\varphi) = (0.01 + 5000k, 0.01 + 200l)$ ($k, l = 0,1,2,...,20$).